\documentclass[authordate1-4]{article}

\usepackage{bm}
\usepackage{amsmath}
\usepackage[title]{appendix}
\usepackage{amsthm}
\usepackage{natbib}
\usepackage{amsfonts}
\usepackage{graphicx}
\usepackage{afterpage}
\usepackage{url}
\usepackage{xcolor}
\usepackage[caption=false]{subfig}%
\usepackage[top=0.95in, bottom=1in, left=1.2in, right=1.2in]{geometry}
\newcommand{\bfS}{\mathbf{S}}
\newcommand{\bfH}{\mathbf{H}}
\newcommand{\bfR}{\mathbf{R}}
\newcommand{\bfC}{\mathbf{C}}

\newcommand{\bfB}{\mathbf{B}}

\newcommand{\bfx}{\mathbf{x}}
\newcommand{\bfy}{\mathbf{y}}

\newcommand{\Rp}{\bfR_{RR}}
\newcommand{\Cp}{\bfC_{RR}}
\newcommand{\Rm}{\bfR_{ME}}

\newcommand{\bfV}{\mathbf{V}}
\newcommand{\bfGam}{\boldsymbol{\Gamma}}
\newcommand{\bfLam}{\boldsymbol{\Lambda}}

\newcommand{\bfSig}{\boldsymbol{\Sigma}}
\newcommand{\sigp}{\boldsymbol{\Sigma}_{RR}}
\newcommand{\sigm}{\boldsymbol{\Sigma}_{ME}}
\newcommand{\halfblankline}{\quad\vspace{+0.5\baselineskip}\pagebreak[3]}
\newcommand{\reviews}{\color{black}}

\newtheorem{mydef}{Definition}
\newtheorem{thm}{Theorem}
\newtheorem{cor}{Corollary}

\title{Improving the condition number of estimated covariance matrices}
\author{ Jemima M. Tabeart  \thanks{School of Mathematical, Physical and Computational Sciences, University of Reading, UK, NERC National Centre for Earth Observation, ({jemima.tabeart@pgr.reading.ac.uk})},  \and Sarah L. Dance \thanks{School of Mathematical, Physical and Computational Sciences, University of Reading, UK} \and Amos S. Lawless \footnotemark[1], \and Nancy K. Nichols \footnotemark[1], \and Joanne A. Waller \footnotemark[2]}

\begin{document}
%\graphicspath{{/home/jemima/Dropbox/UoR/PhD/Papers/Reconditioning/}}
\maketitle

\flushleft
\begin{abstract} %150-250 words
{ High dimensional error covariance matrices and their inverses are used to weight the contribution of observation and background information in data assimilation procedures. As observation error covariance matrices are often obtained by sampling methods, estimates are often degenerate or ill-conditioned, making it impossible to invert an observation error covariance matrix without the use of techniques to reduce its condition number.}
In this paper we present new theory for two existing methods that can be used to `recondition' any covariance matrix: ridge regression, and the minimum eigenvalue method.  { We compare these methods with multiplicative variance inflation, which cannot alter the condition number of a matrix, but is often used to account for neglected correlation information.}
We investigate the impact of reconditioning on variances and correlations of a general covariance matrix in both a theoretical and practical setting. Improved theoretical understanding provides guidance to users regarding method selection, and choice of target condition number.
The new theory shows that, for the same target condition number, both methods increase variances compared to the original matrix, with larger increases for ridge regression than the minimum eigenvalue method. 
We prove that the ridge regression method strictly decreases the absolute value of off-diagonal correlations. { Theoretical comparison of the impact of reconditioning and multiplicative variance inflation on the data assimilation objective function shows that variance inflation alters information across all scales uniformly, whereas reconditioning has a larger effect on scales corresponding to smaller eigenvalues.}
We then consider two examples: a general correlation function, and an observation error covariance matrix arising from interchannel correlations. The minimum eigenvalue method results in smaller overall changes to the correlation matrix than ridge regression, but can increase off-diagonal correlations. 
{\reviews  Data assimilation experiments reveal that reconditioning corrects spurious noise in the analysis but underestimates the true signal compared to multiplicative variance inflation.  }
%We find cases where standard deviations doubled in size, which will have implications for the analysis of a data assimilation problem.

\end{abstract}

\section{Introduction}

The estimation of covariance matrices for large dimensional problems is of growing interest \citep{Pourahmadi13}, particularly for the field of numerical weather prediction (NWP) \citep{bormann16,weston14} where error covariance estimates are used as weighting matrices in data assimilation problems, {\reviews e.g. \citet{Daley92,GHIL1989171,GHIL1991141}}.  At operational NWP centres there are typically $\mathcal{O}(10^7)$ measurements every 6 hours \citep{Bannister17}, meaning that observation error covariance matrices  are extremely high-dimensional. 
In nonlinear least squares problems arising in variational data assimilation, the inverse of correlation matrices are used, meaning that well-conditioned matrices are vital for practical applications \citep{Bannister17}. 
{This is true in both the unpreconditioned and preconditioned variational data assimilation problem using the control variable transform, as the inverse of the observation error covariance matrix appears in both formulations.} {\reviews The convergence of the data assimilation problem can be poor if either the background or observation variance is small; however, }
the condition number and eigenvalues of {\reviews background and} observation error covariance matrices have also been shown to be important for convergence {in both the unpreconditioned and preconditioned case} in  \citet{haben11b,haben11a,haben11c,tabeart17a}. 
Furthermore, the conditioning and solution of the data assimilation system can be affected by complex interactions between the background and observation error covariance matrices and the observation operator \citep{tabeart17a,johnson05}.
The condition number of a matrix, $\textbf{A}$, provides a measure of the sensitivity of the solution $\bfx$ of the system $\textbf{A}\bfx=\mathbf{b}$ to perturbations in $\mathbf{b}$. {\reviews The need for well-conditioned background and observation error covariance matrices} motivates the use of `reconditioning' methods, which are used to reduce the condition number of a given matrix. \halfblankline

{ In NWP applications, observation error covariance matrices are  often constructed from a limited number of samples {\reviews \cite{cordoba16,waller16b,waller16c}}. 
This can lead to problems with sampling error, which manifest in sample covariance matrices, or other covariance matrix estimates, that are very ill-conditioned or can fail to satisfy required properties of covariance matrices (such as symmetry and positive semi-definiteness) \citep{Higham16,Ledoit04}. %In NWP the diagnostic of \citet{desroziers05} is commonly used to estimate observation error covariances by calculating observation-minus-background and observation-minus-analysis residuals \citep{cordoba16,waller16b,waller16c}).
%In practice diagnostic methods can fail to account for some sources of error, resulting in poor estimates of error covariances. 
In some situations it may be possible to determine which properties of the covariance matrix are well estimated. One such instance is presented in \citet{Skoien2006}, which considers how well we can expect the mean, variance and correlation lengthscale of a sample correlation to represent the true correlation matrix depending on different properties of the measured domain (e.g. sample spacing, area measured by each observation). However, this applies only to direct estimation of correlations and will not apply to  diagnostic methods, e.g. \citet{desroziers05}, where transformed samples are used {\reviews and covariance estimates may be poor}. {\reviews We note that} in this paper, we assume that the estimated covariance matrices {\reviews used in our experiments} represent the desired correlation information matrix well and that differences are due to noise rather than neglected sources of uncertainty. {\reviews This may not be the case for practical situations, where reconditioning may need to be performed in conjunction with other techniques to compensate for the underestimation of some sources of error.} \halfblankline

%\citet{waller16,waller17} provide guidance on how the diagnosed matrices are likely to differ from the true statistics in this case, given an understanding of the assimilation system used to produce the residual samples.\halfblankline

 %Interaction between mis-specified  variances and lengthscales for both background and observation error covariances can make it difficult to quantify the impact on the estimated observation error covariance matrices.  Complex interactions between background and observation error covariance matrices and the observation operator are investigated explicitly in \citet{johnson05} for the incremental 4D-Var problem.  In \citet{tabeart17a}, these interactions are shown also to affect the conditioning of the assimilation system.  . {\reviews This may not be the case for practical situations, where reconditioning may need to be performed in conjunction with other techniques to compensate for the underestimation of some sources of error. }}\halfblankline

 { % It is difficult to avoid sampling noise as most estimates of error statistics are obtained using statistical samples.
 	 Depending on the application, a variety of methods have been used to combat the problem of rank deficiency of sample covariance matrices. In the case of spatially correlated errors it may be possible to fit a smooth correlation function or operator to the sample covariance matrix as was done in \citet{simonin18,Guillet19} respectively. Another approach is to retain only the first $k$ leading eigenvectors of the estimated correlation matrix and to add a diagonal matrix to ensure the resulting covariance matrix has full rank \citep{Michel18,stewart13}. However, this has been shown to introduce noise at large scales for spatial correlations and may be expensive in terms of memory and computational efficiency \citep{Michel18}. Although localisation can be used to remove spurious correlations, and can also be used to increase the rank of a degenerate correlation matrix \citep{Hamill01}, it struggles to reduce the condition number of a matrix without destroying off-diagonal correlation information \citep{Smith17}. 
	A further way to increase the rank of a matrix is by considering a subset of columns of the original matrix that are linearly
	independent. This corresponds to using a subset of observations, which is contrary to a key
	motivation for using correlated observation error statistics: the ability to include a larger number of observations in the assimilation system \citep{janjic17}.  Finally, the use of transformed observations imay result in independent {\reviews observation} errors \citep{migliorini12,prates16}{\reviews; however}, %, off-diagonal correlation information is moved into the transformed observation and transformed observation operator. This 
	%means that
	 problems with conditioning will manifest in other components of the data assimilation algorithm, typically the observation operator. Therefore, although other techniques to tackle the problem of ill-conditioning exist, they each have limitations. This suggests that for many applications the use of reconditioning methods, which we will show are inexpensive to implement and are not limited to spatial correlations, may be beneficial.}\halfblankline

 %Theoretical results presented in \citet{haben11a,haben11b,haben11c,tabeart17a} demonstrate that the eigenvalues of both background and observation error covariance matrices are important for
%determining the conditioning of the variational objective function, so it is expected that
%reconditioning of both error covariance matrices would improve the conditioning of the overall
%problem. %However, as the covariance matrix is altered by reconditioning methods, we do not wish
%to use such methods unless they are necessary (as is the case for the observation error covariance matrix).
%\cite{campbell16} indicate that application of reconditioning methods to
%observation error covariance matrices also improves convergence of the dual variational data
%assimilation formulation. 
{\reviews We note that small eigenvalues of the observation error covariance matrix are not the only reason for slow convergence: if observation standard deviations are small, the observation error covariance matrix may be well-conditioned, but convergence of the minimisation problem is likely to be poor \citep{haben11c,tabeart17a}. In this case reconditioning may not improve convergence and performance of the data assimilation routine. } \halfblankline

 Two methods in particular, referred to in this work as the minimum eigenvalue method and ridge regression, are commonly used { at NWP centres}.
  Both methods are { used by} \citet{weston11}, where they are tested numerically. Additionally in \citet{campbell16} a comparison between these methods is made experimentally {\reviews and it is shown that reconditioning improves convergence of a dual} four-dimensional variational assimilation system.  However, up to now there has been minimal theoretical investigation into the effects of these methods on the covariance matrices. {\reviews In this paper we develop theory that shows how variances and correlations are altered by the application of reconditioning methods to a  covariance matrix.}\halfblankline% Better understanding of how these methods change the covariance matrix will allow users to make informed decisions about which method is most appropriate for their application. {\reviews Furthermore, increased theoretical understanding will allow either reconditioning method to be applied in a more rigorous manner.}\halfblankline 
 
%Another benefit of increased theoretical understanding is the opportunity to implement either reconditioning method in a more rigorous manner. Both methods alter the eigenvalues of the original covariance matrix in order to obtain an approximation to the covariance matrix. This change in eigenvalues is determined by a user-specified condition number.
% Currently this condition number is chosen in an ad hoc manner and investigation into the underlying mathematical theory of both methods may provide insight into procedures for choosing the condition number. \halfblankline 

{\reviews Typically reconditioning is applied to improve convergence of a data assimilation system by reducing the condition number of a matrix. However, the convergence of a data assimilation system can also be improved using} 
{multiplicative variance inflation,  a commonly used method at NWP centres such as ECMWF \citep{Liu03,McNally06,bormann15,bormann16} to account for neglected error correlations {\reviews or to address deficiencies in the estimated error statistics} by increasing the uncertainty in observations. It is not a method of reconditioning {\reviews when a constant inflation factor is used}, as it cannot change the condition number of a covariance matrix. {\reviews In practice multiplicative variance inflation is often combined with other techniques, such as neglecting off-diagonal error correlations, which do alter the conditioning of the observation error covariance matrix}.\halfblankline
	%	{\reviews Multiplicative variance inflation can improve convergence of a data assimilation problem, but the motivation for applying this method is primarily to mitigate for underestimation of error statistics, for example when using uncorrelated observation error matrices. In these cases, reconditioning methods may require very low target condition numbers, or be inappropriate to use alone. In particular, reconditioning works well when  sample estimates provide a good estimate of the true error statistics; it  can be hard for users to assess whether this is the case for their system. } \halfblankline

	 {\reviews Although it is not a reconditioning technique,} in \citet{bormann15} multiplicative variance inflation was found  to yield faster convergence of a data assimilation procedure than either the ridge regression or minimum eigenvalue methods of reconditioning. %to improve the convergence rate of the assimilation system. %
	This finding is likely to be system-dependent; the original diagnosed error covariance matrix in the ECMWF system has a smaller condition number than the corresponding matrix for the same instrument in the Met Office system \citep{weston14}. 
%	{\reviews SOMETHING ABOUT }
	Additionally, in the ECMWF system the use of reconditioning methods only results in small improvements to convergence, and there is little difference in convergence speed for the two methods. This contrasts with the findings of \citet{weston11,weston14,campbell16} where differences in convergence speed when using each method of reconditioning were found to be large. Therefore, it is likely that the 
	importance of reducing the condition number of the observation error covariance matrix compared to inflating variances will be sensitive to the data assimilation system of interest. 
 {\reviews Aspects of the data assimilation system that may be important in determining the level of this sensitivity include: the choice of preconditioning and minimisation scheme \citep{bormann15}, quality of the covariance estimate, interaction between background and estimated observation error covariance matrices within the data assimilation system \citep{fowler17,tabeart17a}, the use of thinning and different observation networks.} 
 	 We also note that \citet{stewart08,stewart09,stewart13} consider changes to the information content and analysis accuracy corresponding to different approximations to a correlated observation error covariance matrix (including an inflated diagonal matrix). \citet{stewart13,healy05} also provide evidence in idealized cases to show that inclusion of even approximate correlation structure gives significant benefit over diagonal matrix approximations, including when variance inflation is used.  }\halfblankline

In this work we investigate the minimum eigenvalue and ridge regression methods of reconditioning { as well as multiplicative variance inflation}, and analyse their impact on the covariance matrix. We compare both methods theoretically for the first time, by considering the impact of reconditioning on the correlations and variances of the covariance matrix. { We also study how each method alters the objective function when applied to the observation error covariance matrix.} Other methods of reconditioning, including thresholding \citep{Bickel08} and localisation \citep{Horn91,Menetrier15,Smith17} %linear shrinkage \citep{Ledoit04} %and regularisation methods such as the Lasso penalty technique \citep{Pourahmadi13}, 
have been discussed from a theoretical perspective in the literature but will not be included in this work. In Section \ref{sec:Methods} we describe the methods more formally than in previous literature before developing new related theory in detail in Section \ref{sec:RecondTheory}. We show that the ridge regression method increases the variances and decreases the correlations for a general covariance matrix and the minimum eigenvalue method increases variances. We prove that the increases to the variance are bigger for the ridge regression method than the minimum eigenvalue method for any covariance matrix. { We show that both methods of reconditioning reduce the weight on observation information in the objective function in a scale dependent way, with the largest reductions in weight corresponding to the smallest eigenvalues of the original observation error covariance matrix. In contrast, multiplicative variance inflation {\reviews using a constant inflation factor} reduces the weight on observation information by a constant amount for all scales. } In Section \ref{sec:NumExpt} the methods are illustrated via numerical experiments for two types of covariance structures. One of these is a simple general correlation function, and one is an interchannel covariance arising from a satellite based instrument with observations used in NWP. { We provide physical interpretation of how each method alters the covariance matrix, and use this to provide guidance on which method of reconditioning is most appropriate for a given application.} {\reviews We present an illustration of how all three methods alter the analysis of a data assimilation problem, and relate this to the theoretical conclusions concerning the objective function. }We finally present our conclusions in Section \ref{sec:Conclusions}. The methods are very general and, although their initial application was to observation error covariances arising from numerical weather prediction, the results presented here apply to any sampled covariance matrix, such as those arising in %{\reviews oceanography \citep{RobertsJones16}, ecology \citep{COOPER2018199,PinningtonEwanM.2016Itro,PinningtonEwanM.2017Uteo}}, 
finance \citep{Higham02,Qi2010} and neuroscience \citep{Nakamura15,SchiffStevenJ2011NCET}.

{\reviews\section{Covariance matrix modification methods} \label{sec:Methods}}

We begin by defining the condition number.
{\reviews The condition number provides a measure of how sensitive the solution $\mathbf{x}$ of a linear equation $\mathbf{A}\bfx = \mathbf{b}$ is to perturbations in the data $\mathbf{b}$. A `well-conditioned problem' will result in small perturbations to the solution with small changes to  $\mathbf{b}$, whereas for an `ill-conditioned problem', small perturbations to  $\mathbf{b}$ can result in large changes to the solution.} , Noting that all covariance matrices are positive semi-definite by definition, we distinguish between the two cases of strictly positive definite covariance matrices and covariance matrices with zero minimum eigenvalue. Symmetric positive definite matrices admit a definition for the condition number in terms of their maximum and minimum eigenvalues. For the remainder of the work, we define the eigenvalues of a symmetric positive semi-definite matrix $\bfS\in \mathbb{R}^{d\times d}$ via:
\begin{equation}\label{eq:evalordering}
 \lambda_{1}(\bfS) \ge \dots\ \ge \lambda_d(\bfS)\ge 0.
\end{equation} 

\begin{thm} \label{def:condno}
	If $\bfS \in \mathbb{R}^{d \times d}$ is a symmetric positive definite matrix with eigenvalues defined as in \ref{eq:evalordering} we can write the condition number in the $L_2$ norm as 
%\begin{equation}
$	\kappa(\bfS)= \frac{\lambda_{1}(\bfS)}{\lambda_{d}(\bfS)}.$
%	\end{equation}
%	where $\lambda_{1}(\bfS)$ and $\lambda_{N}(\bfS)$ correspond to the largest and smallest eigenvalues of $\bfS$ respectively. 
\end{thm}
\begin{proof}
	See \citep[Sec. 2.7.2]{golub96}.
\end{proof}
For a singular covariance matrix, $\bfS$, the convention is to take $\kappa(\bfS) = \infty$ \citep[Sec. 12]{TrefethenLloydN.LloydNicholas1997Nla}. We also note that real symmetric matrices admit orthogonal eigenvectors which can be normalised to produce a set of orthonormal eigenvectors. \halfblankline

Let $\bfR \in \mathbb{R}^{d\times d}$ be a positive semi-definite covariance matrix with condition number $\kappa(\bfR) = \kappa$. We wish to recondition $\bfR$ to obtain a covariance matrix with condition number $\kappa_{max}$
 \begin{equation}\label{eq:kappadefn}
 1 \le\kappa_{max}<\kappa,
 \end{equation} where the value of $\kappa_{max}$ is chosen by the user.  We denote the eigendecomposition of $\bfR$ by 

\begin{equation}\label{eq:eigendecomp}
\bfR=\bfV_{\bfR}\bfLam\bfV_{\bfR}^T
\end{equation}
 where $\bfLam\in\mathbb{R}^{d\times d}$ is the diagonal matrix of eigenvalues of $\bfR$ and $\bfV_{\bfR}\in \mathbb{R}^{d\times d}$ is a corresponding matrix of orthonormal eigenvectors.\halfblankline

In addition to considering how the covariance matrix itself changes with reconditioning, it is also of interest to consider how the related correlations and standard deviations are altered. We decompose $\bfR$ as 
$\bfR =  \bfSig \mathbf{C}  \bfSig,$
where $\bfC$ is a correlation matrix, and $ \bfSig$ is a non-singular diagonal matrix of standard deviations. We calculate $\bfC$ and $ \bfSig$ via:
\begin{equation} \label{eq:CovtoVar}
\bfSig(i,i) = \sqrt{\bfR(i,i)}, \quad
\bfC(i,j) = \frac{\bfR(i,j)}{\sqrt{\bfR(i,i)}\sqrt{\bfR(j,j)}}  .
\end{equation}

We now introduce the ridge regression method and the minimum eigenvalue method; the two methods of reconditioning that will be discussed in this work. { We then define  multiplicative variance inflation. This last method is not a method of reconditioning, but will be used for comparison purposes with the ridge regression and minimum eigenvalues methods.}
\subsection{Ridge regression method:} \label{sec:RRalg}
The ridge regression method (RR) adds a scalar multiple of the identity to $\bfR$ to obtain the reconditioned matrix $\bfR_{RR}$.  The scalar $\delta$ is set using the following method.
\begin{itemize} %Layout of this - do we want it all on one line? Could we italicise and centre?
	\item Define $\delta = \frac{\lambda_{1}(\bfR) - \lambda_{d}(\bfR)\kappa_{max}}{\kappa_{max}-1} $. \linebreak
	\item Set $\Rp = \bfR+\delta\mathbf{I}$
\end{itemize}%\halfblankline

We note that this choice of $\delta$ yields $\kappa(\Rp)=\kappa_{max}$. \halfblankline

{ In the literature \citep{Hoerl70,Ledoit04}, `ridge regression' is a method used to regularise least squares problems. In this context, ridge regression can be shown to be equivalent to Tikhonov regularisation \citep{Hansen98}. However, in this paper we apply ridge regression as a reconditioning method directly to a covariance matrix. For observation error covariance matrices, the reconditioned matrix is then inverted prior to its use as a weighting matrix in the data assimilation objective function. As we are only applying the reconditioning to a single component matrix in the variational formulation, the implementation of the ridge regression method used in this paper is not equivalent to Tikhonov regularisation applied to the variational data assimilation problem \citep{budd11,Moodey_2013}. This is shown in Section	\ref{sec:ObjectiveFn} where we consider how applying ridge regression to the observation error covariance matrix affects the variational data assimilation objective function.} The ridge regression method is used at the Met Office \citep{weston14}.
\subsection{Minimum eigenvalue method:} \label{sec:MinEig}

The minimum eigenvalue method (ME) fixes a threshold, $T$, below which all eigenvalues of the reconditioned matrix, $\bfR_{ME}$, are set equal to the threshold value. The value of the threshold is set using the following method.

\begin{itemize}

	\item Set $\lambda_{1}(\Rm) = \lambda_1(\bfR)$\linebreak
	\item Define $T={\lambda_{1}(\bfR)}/{\kappa_{max}}>\lambda_d(\bfR)$, where $\kappa_{max}$ is defined in \eqref{eq:kappadefn}.\linebreak
	\item Set the remaining eigenvalues of $\Rm$ via	\begin{equation} \label{eq:mineigupdate}
	\lambda_{k}(\Rm)= 
	\begin{cases}
	\lambda_k(\bfR) & \text{if } \lambda_k(\bfR)>T\\
	T & \text{if } \lambda_k(\bfR) \le T \\
	\end{cases} .
	\end{equation}
	\item Construct the reconditioned matrix via $\Rm = \mathbf{V}_{\bfR}\bfLam_{ME}\mathbf{V}_{\bfR}^T$, where $\bfLam_{ME}(i,i) = \lambda_i(\Rm)$. 

\end{itemize}%\halfblankline

This yields $\kappa(\Rm) = \kappa_{max}$.
The updated matrix of eigenvalues can be written as $\bfLam_{ME} = \bfLam + \bfGam$, the sum of the original matrix of eigenvalues and $\bfGam$, a low-rank diagonal matrix update with entries $\bfGam(k,k) = \max\{T - \lambda_k, 0\}$. Using \eqref{eq:eigendecomp} the reconditioned $\Rm $ can then be written as:
\begin{equation}\label{eq:MEdecompalt}
\Rm = \bfV_{\bfR}(\bfLam + \bfGam)\bfV^T_{\bfR} = \bfR+ \bfV_{\bfR} \bfGam \bfV^T_{\bfR}.
\end{equation}  

Under the condition that $\kappa_{max}>d-l+1$, where $l$ is the index such that $\lambda_l\le T <\lambda_{l-1}$, %that $\kappa_{max}\ge d$ 
the minimum eigenvalue method is equivalent to minimising 
the difference $\bfR-\Rm \in \mathbb{R}^{d\times d}$ with respect to the Ky Fan 1-$d$ norm (The proof is provided in Appendix \ref{sec:Append}). {The Ky Fan $p-k$ norm (also referred to as the trace norm) is defined in \citet{KyFan59,Horn91}, and is used in \citet{takana14} to find the closest positive definite matrix with condition number smaller than a given constant.} A variant of the minimum eigenvalue method is applied to observation error covariance matrices at the European Centre for Medium-Range Weather Forecasts (ECMWF) \citep{bormann16}.

{\subsection{Multiplicative variance inflation}\label{sec:MultInfl}
 
 Multiplicative variance inflation (MVI) is a method that increases the variances corresponding to a covariance matrix. Its primary use is to account for neglected error correlation information, particularly in the case where diagonal covariance matrices are being used even though non-zero correlations exist in practice. However, this method can also be applied to non-diagonal covariance matrices. 
 
 \begin{mydef}\label{def:MVI}
 	Let $\alpha>0$ be a given variance inflation factor and $$\bfR=\bfSig\mathbf{C}\bfSig$$ be the estimated covariance matrix. Then multiplicative variance inflation is defined by \begin{equation}
 	\bfSig_{MVI} = \alpha\bfSig.
 	\end{equation}
 	This is equivalent to multiplying the estimated covariance matrix by a constant. The updated covariance matrix is given by 
 	\begin{equation}
 	\bfR_{MVI} =\left(\alpha\bfSig\right)\bfC\left(\alpha\bfSig\right)= \alpha^2\bfSig\bfC\bfSig=\alpha^2\bfR.
 	\end{equation}
 	The estimated covariance matrix is therefore multiplied by the square of the inflation constant.
 	We note that the correlation matrix, $\bfC$, is unchanged by application of multiplicative variance inflation. 
 \end{mydef}

Multiplicative variance inflation is used at NWP centres including ECMWF \citep{bormann16} {\reviews to counteract deficiencies in estimated error statistics, such as underestimated or neglected sources of error. Inflation factors are tuned to achieve improved analysis or forecast performance, and are hence strongly dependent on the specific data assimilation system of interest. Aspects of the system that might influence the choice of inflation factor include observation type, known limitations of the covariance estimate, and observation sampling or thinning.} \halfblankline

Although variance inflation is not a method of reconditioning, as it is not able to alter the condition number of a covariance matrix, we include it in this paper for comparison purposes. This means that variance inflation can only be used in the case that the estimated matrix can be inverted directly, i.e. is full rank. Multiplicative variance inflation could also refer to the case where the constant inflation factor is replaced with a diagonal matrix of inflation factors. In this case the condition number of the altered covariance matrix would change. {\reviews An example of a study where multiple inflation factors are used is given by \citet{Heilliette15}, where the meteorological variable to which an observation is sensitive determines the choice of inflation factor. }  However, this is not commonly used in practice, and will not be considered in this paper.}

\section{{ Theoretical considerations}}\label{sec:RecondTheory}

In this section we develop new theory for each method. We are particularly interested in the changes made to $\bfC
$ and $ \bfSig$ for each case. 
Increased understanding of the effect of each method may allow users to adapt or extend these methods, or determine which is the better choice for practical applications. \halfblankline

We now introduce an assumption that will be used in the theory that follows.\halfblankline

\textbf{Main Assumption:} Let $\bfR\in \mathbb{R}^{d \times d}$ be a symmetric positive semi-definite matrix with $\lambda_1(\bfR)>\lambda_d(\bfR)$.\halfblankline

We remark that any symmetric, positive semi-definite matrix with $\lambda_1=\lambda_d$ is a scalar multiple of the identity, and cannot be reconditioned since it is already at its minimum possible value of unity. Hence in what follows, we will consider only matrices $\bfR$ that satisfy the Main Assumption.

\subsection{Ridge Regression Method} \label{sec:RidgeRegression}

We begin by discussing the theory of RR. In particular we prove that applying this method { for any positive scalar, $\delta$,} results in a decreased condition number for any choice of $\bfR$.
\begin{thm}\label{lem:ridgeregincreasecond}
	
	Under the conditions of the Main Assumption, adding { any} positive increment to the diagonal elements of $\bfR$ decreases its condition number.
\end{thm}
\begin{proof}

	We recall that $\Rp=\bfR+\delta\mathbf{I}$. The condition number of $\Rp$ is given by
	\begin{equation} \label{eq:kappaRp}
	\kappa(\Rp) = \frac{\lambda_{1}(\Rp)}{\lambda_{d}(\Rp)} = \frac{\lambda_{1}(\bfR)+\delta}{\lambda_{d}(\bfR)+\delta} .
	\end{equation}
	It is straightforward to show that for any $\delta>0$, $\kappa(\Rp) < \kappa(\bfR)$, completing the proof.  

\end{proof}

We now consider how application of RR affects the correlation matrix $\bfC$ and the diagonal matrix of standard deviations $ \bfSig$.

\begin{thm} \label{lem:updateCandSigRR}

	Under the conditions of the Main Assumption, the ridge regression method updates the standard deviation matrix $\sigp$, and correlation matrix $\Cp$ of $\Rm$ via
	\begin{equation}\label{eq:RRSigupdate}
\sigp = ( \bfSig^2+\delta \mathbf{I}_d  )^{1/2}, \quad
\Cp = \sigp^{-1}\bfR\sigp^{-1} + \delta\sigp^{-2}.
	\end{equation} 

\end{thm}
\begin{proof}
	Using \eqref{eq:CovtoVar}, $ \bfSig(i,i) = (\bfR(i,i))^{1/2}$. Substituting this into the expression for $\Rp$ yields:
	\begin{equation}\label{eq:Sigupdate}
	\sigp(i,i) = (\Rp(i,i))^{1/2} = (\bfR(i,i) + \delta)^{1/2} = ( \bfSig(i,i)^2+ \delta)^{1/2}.
	\end{equation}

Considering the components of $\Cp$ and the decomposition of $\sigp$ given by \eqref{eq:CovtoVar}: 
	\begin{equation}
\sigp\Cp\sigp = \bfR + \delta \mathbf{I}_d, \quad
	\Cp = \sigp^{-1}\bfR\sigp^{-1}+\delta\sigp^{-2}
	\end{equation}
as required.  
\end{proof}

Theorem \ref{lem:updateCandSigRR} shows how we can apply RR to our system by updating $\bfC$ and $ \bfSig$ rather than $\bfR$.
We observe, from \eqref{eq:RRSigupdate}, that applying RR leads to a { constant increase to variances} for all variables. {However, the inflation to standard deviations is additive, rather than the multiplicative inflation that occurs for multiplicative variance inflation. }We now show that RR also reduces all non-diagonal entries of the correlation matrix.   

\begin{cor}
	Under the conditions of the Main Assumption, for $i \ne j$, $|\Cp(i,j)| < |\bfC(i,j)|$.
\end{cor}
\begin{proof}
	Writing the update equation for $\bfC$, given by \eqref{eq:RRSigupdate}, in terms of the variance and correlations of $\bfR$ yields:
	\begin{equation}\label{eq:eqCtemp}
	\Cp = \sigp^{-1} \bfSig \bfC  \bfSig\sigp^{-1} + \delta\sigp^{-2}.
	\end{equation}
	We consider $\Cp(i,j)$ for $i \ne j$. As $\sigp$ and $\bfSig$ are diagonal matrices, we obtain
	\begin{equation}\label{eq:Ctempup}
	\Cp(i,j)
	= \sigp^{-1}(i,i) \bfSig(i,i) \bfC(i,j)  \bfSig(j,j)\sigp^{-1}(j,j).
	\end{equation}
	From the update equation \eqref{eq:RRSigupdate},  $\sigp(i,i) >  \bfSig(i,i)$ for any choice of $i$.
	This means that $\sigp^{-1}(i,i)  \bfSig(i,i) <1$ for any choice of $i$. Using this in \eqref{eq:Ctempup} yields that for all values of $i,j$ with $i\ne j$, $|\Cp(i,j)| < |\bfC(i,j)|$ as required.\halfblankline
	
	For $i=j$, it follows from \eqref{eq:eqCtemp} that $\Cp(i,i)=1$ for all values of $i$.   
\end{proof}

\subsection{Minimum Eigenvalue Method}\label{sec:MEtheory}
We now discuss the theory of ME as introduced in Section \ref{sec:MinEig}.
Using the alternative decomposition of $\Rm$ given by \eqref{eq:MEdecompalt} enables us to update directly the standard deviations for this method.

\begin{thm} \label{lem:changetoSigME}
	Under the conditions of the Main Assumption, the minimum eigenvalue method updates the standard deviations, $\sigm$, of $\bfR$ via
	\begin{equation}\label{eq:MESigupdate}
	\sigm(i,i) = \left(\bfR(i,i) + \sum_{k=1}^{d}\bfV_R(i,k)^2\bfGam(k,k) \right)^{1/2} .
	\end{equation}
	This can be bounded by 
	\begin{equation}\label{eq:MESigRecondbound}
	\bfSig(i,i) \le \sigm(i,i) \le \left(\bfSig(i,i)^2 + T-\lambda_{d}(\bfR) \right)^{1/2}.
	\end{equation}
\end{thm}
\begin{proof}
	\begin{align}
	\sigm(i,i)
	& = \left(\bfR(i,i) + \left(\bfV_R\bfGam\bfV_R^T\right)(i,i)\right)^{1/2} \label{eq:RdecompME}\\ 
	& = \left(\bfR(i,i) + \sum_{k=1}^{d}\bfV_R(i,k)^2\bfGam(k,k) \right)^{1/2} 
	\end{align}
	Noting that $\bfGam(k,k)\ge 0$ for all values of k, we bound the second term in this expression by
	\begin{align}
	0 \le \sum_{k=1}^{d}\bfV_R(i,k)^2\bfGam(k,k) &\le \max_k\{\bfGam(k,k)\}\sum_{k=1}^{d}\bfV_R(i,k)^2 \\
	& \le\left(T-\lambda_d(\bfR)\right) \sum_{k=1}^{d}\bfV_R(i,k)^2 
	 \le T - \lambda_d(\bfR)
	\end{align} 
	This inequality follows from the orthonormality of $\bfV_{\bfR}$, and by the fact that $T>\lambda_d(\bfR)$ by definition.
	 
\end{proof}
Due to the way the spectrum of $\bfR$ is altered by ME it is not evident how correlation entries are altered in general for this method of reconditioning.  

{ \subsection{Multiplicative variance inflation}
We now discuss theory of MVI that was introduced in Section~\ref{sec:MultInfl}. 
We prove that MVI is not a method of reconditioning, as it does not change the condition number of a covariance matrix.
\begin{thm}\label{lem:VarInflRank}
	Multiplicative variance inflation {\reviews with a constant inflation parameter} cannot change the condition number or rank of a matrix.
\end{thm}
\begin{proof}
	Let $\alpha^2>0$ be our multiplicative inflation constant such that $
	\bfR_{MVI} = \alpha^2\bfR$.
	The eigenvalues of $\bfR_{MVI}$ are given by $\alpha^2\lambda_{1}, \alpha^2\lambda_2, \dots, \alpha^2\lambda_{d}$.\halfblankline
	
	If $\bfR$ is rank-deficient, then $\lambda_{min} (\bfR_{MVI}) =\alpha^2\lambda_{d} = 0$ and hence $\bfR_{MVI}$ is also rank deficient.
	If $\bfR$ is full rank then we can compute the condition number of $\bfR_{MI}$ as the ratio of its eigenvalues, which yields
	\begin{equation}
	\kappa(\bfR_{MVI}) = \frac{\alpha^2\lambda_{1}}{\alpha^2\lambda_{d}} 	= \kappa(\bfR)
	%&	= \frac{\lambda_{1}}{\lambda_{d}}\\
	\end{equation}
	Hence the condition number and rank of $\bfR$ are unchanged by multiplicative inflation.
\end{proof}}
\subsection{Comparing ridge regression and minimum eigenvalue methods}
Both RR and ME change $\bfR$ by altering its eigenvalues. In order to compare the two methods, we can consider their effect on the standard deviations. We recall from Sections \ref{sec:RidgeRegression} and \ref{sec:MEtheory} that RR increases standard deviations by a constant and the changes to standard deviations by ME can be bounded above and below by a constant.

\begin{cor}\label{lem:RRvsMEvariance}
	Under the conditions of the Main Assumption, for a fixed value of $\kappa_{max}<\kappa$, $\sigm(i,i) < \sigp(i,i)$ for all values of $i$.
\end{cor}
\begin{proof}
	From Theorems \ref{lem:updateCandSigRR} and \ref{lem:changetoSigME} the updated standard deviation values are given by
	\begin{equation}
	\sigp = \left( \bfSig^2+\delta \mathbf{I}_d \right)^{1/2}
\quad \text{and} \quad \quad
	\sigm(i,i)\le  \left(\bfSig(i,i)^2 + T-\lambda_{d}(\bfR) \right)^{1/2}.
	\end{equation}

From the definitions of $\delta$ and $T$ we obtain that 
\begin{equation}
\delta = \frac{\lambda_1(\bfR)-\lambda_d(\bfR)\kappa_{max}}{\kappa_{max}-1} >\frac{\lambda_1(\bfR)- \lambda_d(\bfR)\kappa_{max}}{\kappa_{max}}= T-\lambda_d(\bfR).
\end{equation}
	We conclude that the increment to the standard deviations for RR is always larger than the increment for ME.  
\end{proof}

{ \subsection{Comparison of methods of reconditioning and multiplicative variance inflation on the variational data assimilation objective function}\label{sec:ObjectiveFn}
We demonstrate how RR, ME and MVI alter the objective function of the variational data assimilation problem when applied to the observation error covariance matrix. We consider the 3D-Var objective function here for simplicity of notation, although the analysis extends naturally to the 4D-Var case. 
We begin by defining the 3D-Var objective function of the variational data assimilation problem. 
\begin{mydef}\label{def:VarCostFn}
	The objective function of the variational data assimilation problem is given by 
	\begin{equation}\label{eq:CostFn}
	J(\bfx) =\frac{1}{2} (\bfx-\bfx_b)^T\bfB^{-1}(\bfx-\bfx_b)+\frac{1}{2}(\bfy-h[\bfx])^T\bfR^{-1}(\bfy-h[\bfx]) := J_b+J_o
	\end{equation}
	where $\bfx_b\in \mathbb{R}^n$ is the background or prior, $\bfy\in\mathbb{R}^{d}$ is the vector of observations, $h:\bfR^n\rightarrow\bfR^d$ is the observation operator mapping from control vector space to observation space, $\bfB\in \mathbb{R}^{n\times n}$ is the background error covariance matrix, and $\bfR\in \mathbb{R}^{d\times d}$ is the observation error covariance matrix.
	Let $J_o$ denote the observation term in the objective function and $J_b$ denote the background term in the objective function.
\end{mydef}

In order to compare the effect of using each method, they are applied to the observation error covariance matrix in the variational objective function \eqref{eq:CostFn}. We note that analogous results hold if all methods are applied to the background error covariance matrix in the objective function.\halfblankline %$\bfx,\mathbf{x},x,\mathbf{x},\bm{x}$

{\reviews We begin by presenting  the three updated objective functions, and then discuss the similarities and differences for each method together at the end of Section 3.5.} We first consider how applying RR to the observation error covariance matrix alters the variational objective function \eqref{eq:CostFn}. 
\begin{thm}
	By applying RR to the observation error covariance matrix we alter the objective function \eqref{eq:CostFn} as follows:
	%	\begin{equation}
	\begin{equation}\label{eq:RRCostFn}
	J_{RR}(\bfx)%)&= J_b+J_o-(\bfy-h[\bfx])^TV_{\bfR}\bfLam{\delta} V_{\bfR}^T(\bfy-h[\bfx]) \\
	%	&=J(\bfx) - || V_{\bfR}^T(\bfy-h[\bfx])||_{(\bfLam{\delta})^{1/2}}^2\\
	= J(\bfx)-(\bfy-h[\bfx])^TV_{\bfR}\bfLam_{\delta}V_{\bfR}^T(\bfy-h[\bfx])
	\end{equation}
	%	\end{equation}
	where $\bfLam_{\delta}$ is a diagonal matrix with entries given by $(\bfLam_{\delta})_{ii} = \frac{\delta}{\lambda_i(\lambda_i+\delta)}$.
\end{thm}
\begin{proof}
We denote the eigendecomposition of $\bfR$ as in	\eqref{eq:eigendecomp}.
 Applying RR to the observation error covariance matrix, $\bfR$, we obtain
	\begin{equation}
	\bfR_{RR}=\bfV_{\bfR}(\bfLam+\delta\mathbf{I}_p)\bfV_{\bfR}^T.
	\end{equation}
	
	We then calculate the inverse of $\bfR_{RR}$ and express this in terms of $\bfR^{-1}$ and an update term: 
	\begin{align}
	\bfR_{RR}^{-1} &= \bfV_{\bfR}(\bfLam+\delta\mathbf{I}_p)^{-1}\bfV_{\bfR}^T\\
	&= \bfV_{\bfR}\text{Diag}\Big(\frac{1}{\lambda_i}-\frac{\delta}{\lambda_i(\lambda_i+\delta)}\Big)\bfV_{\bfR}^T\\
	&=\bfR^{-1} - \bfV_{\bfR}\text{Diag}\Big(\frac{\delta}{\lambda_i(\lambda_i+\delta)}\Big)\bfV_{\bfR}^T \label{eq:B12}
	\end{align}

	Substituting \eqref{eq:B12} into  \eqref{eq:CostFn}, and defining $\bfLam_{\delta}$ as in the theorem statement we can write the objective function using the reconditioned observation error covariance matrix as \eqref{eq:RRCostFn}.
\end{proof}
 
{ The effect of RR on the objective function differs from the typical application of Tikhonov regularisation to the variational objective function \citep{budd11,Moodey_2013}. In particular, we subtract a term from the original objective function rather than adding one, and the term depends on the eigenvectors of $\bfR$ as well as the innovations (differences between observations and the background field in observation space). {\reviews Writing the updated objective function as in \eqref{eq:RRCostFn} shows that the size of the original objective function \eqref{eq:CostFn} is decreased when RR is used. Specifically, as we discuss later, the contribution of small-scale information to the observation term, $J_o$, is reduced by the application of RR.} }\halfblankline

We now consider how applying ME to the observation error covariance matrix alters the objective function \eqref{eq:CostFn}.
\begin{thm}
	By applying ME to the observation error covariance matrix we alter the objective function \eqref{eq:CostFn} as follows:
	%	\begin{equation}
	\begin{equation}\label{eq:MEcostfn}
	J_{ME}(\bfx)=J(\bfx)-(\bfy-h[\bfx])^TV_{\bfR}\widetilde{\mathbf{\Gamma}} V_{\bfR}^T(\bfy-h[\bfx]), 
	%&=J(\bfx) - || V_{\bfR}^T(\bfy-h[\bfx])||_{\tilde{\Gamma}^{1/2}}^2\\
%	= J(\bfx)-||\tilde{\bfGam}^{1/2}V_{\bfR}^T(\bfy-h[\bfx])||_2^2.
	\end{equation}
	%	\end{equation}
	where 
	\begin{equation}
	\widetilde{\bfGam}(i,i) = \begin{cases}
	0 & if \lambda_i \ge T\\
	\frac{T-\lambda_i}{T\lambda_i} & if \lambda_i <T
	\end{cases}.
	\end{equation}
\end{thm}
\begin{proof}
	We begin by applying ME and decomposing $\bfR_{ME}$ as in \eqref{eq:MEdecompalt}:
	\begin{equation}
	\bfR_{ME}  = \bfV_{\bfR}(\bfLam+\bfGam)\bfV_{\bfR}^T.
	\end{equation}
	Therefore calculating the inverse of the reconditioned matrix yields
	\begin{equation}
	\bfR_{ME} = \bfV_{\bfR}(\bfLam+\bfGam)^{-1}\bfV_{\bfR}^T
	\end{equation}

	As this is full rank we can calculate the inverse of the diagonal matrix $\bfLam+\bfGam$
	\begin{align}
	(\bfGam+\bfLam)^{-1}(i,i) &= \begin{cases}
	\frac{1}{\lambda_i} & if \lambda_i\ge T\\
	\frac{1}{\lambda_i + (T-\lambda_i)} & if \lambda_i <T
	\end{cases}\\
	& = \bfLam^{-1} - \begin{cases}
	0 & if \lambda_i\ge T\\
	\frac{T-\lambda_i}{T\lambda_i} & if \lambda_i<T
	\end{cases}.
	\end{align}
	
	Defining $\tilde{\bfGam}$ as in the theorem statement, and we can write $\bfR_{ME}^{-1}$ as 
	\begin{equation}
	\bfR_{ME}^{-1} =\bfR^{-1} - \bfV_{\bfR}\tilde{\bfGam}\bfV_{\bfR}^T
	\end{equation}
	Substituting this into the definition of the objective function \eqref{eq:CostFn} we obtain the result given in the theorem statement.
\end{proof}
{\reviews As $\tilde{\bfGam}$ is non-zero only for eigenvalues smaller than the threshold, $T$, the final term of the updated objective function \eqref{eq:MEcostfn} reduces the weight on eigenvectors corresponding to those small eigenvalues. As all the entries of $\tilde{\bfGam}$ are non-negative, the size of the observation term in the original objective function \eqref{eq:CostFn} is decreased when ME is used. }\halfblankline

Finally we consider the impact on the objective function of using MVI. We note that this can only be applied in the case that the estimated error covariance matrix is invertible as, by the result of Theorem \ref{lem:VarInflRank}, variance inflation cannot change the rank of a matrix.

\begin{thm}
	In the case that $\bfR$ is invertible, the application of MVI to the observation error covariance matrix alters the objective function \eqref{eq:CostFn} as follows
	\begin{equation}\label{eq:CostFnMVI}
	J_{MVI}(\bfx) = J_b + \frac{1}{\alpha^2}J_o
	\end{equation}
\end{thm}
\begin{proof}
	By Definition \ref{def:MVI}, $\bfR_{MVI} = \alpha^2\bfR$ for inflation parameter $\alpha$. The inverse of $\bfR_{MVI}$ is given by 
	\begin{equation}
	\bfR_{MVI}^{-1} = \frac{1}{\alpha^2}\bfR^{-1}.
	\end{equation}
	Substituting this into \eqref{eq:CostFn} yields the updated objective function given by \eqref{eq:CostFnMVI}.
\end{proof} %the smallest eigenvalues of R are changed the most relatively 
For both reconditioning methods, {\reviews the largest relative changes to the spectrum of $\bfR$ occur for its smallest eigenvalues. } In the case of {\reviews positive} spatial correlations, {\reviews small eigenvalues} are typically sensitive to smaller scales.  For spatial correlations, weights on scales  of the observations {\reviews associated with smaller eigenvalues} are reduced in the variational objective function, increasing the relative sensitivity of analysis to information content from the observations at large scales.\halfblankline% \citet{fowler17} studied the impact of altering the observation error correlation lengthscales on the analysis, finding that the smallest analysis error variances and greatest spread in information occurs when observation and background are accurate at different scales. %This means that in the case that the observations are more sensitive to smaller scales than the background, reconditioning the observation error covariance matrix will decrease the weight on small scale information and hence increase analysis error variances. \halfblankline

We also see that for RR and ME smaller choices of $\kappa_{max}$ yield larger reductions to the weight applied to small scale observation information. For RR, a smaller target condition number results in a larger value of $\delta$, and hence larger diagonal entries of $\bfLam_{\delta}$. For ME, a smaller target condition number yields a larger threshold, $T$, and hence larger diagonal entries of $\tilde{\bfGam}$. This means that the more reconditioning that is applied, the less weight the observations will have in the analysis. This reduction in observation weighting is different for the two methods; RR reduces the weight on all observations, although the relative effect is larger for scales corresponding to the smallest eigenvalues, whereas ME only reduces weight for scales corresponding to eigenvalues smaller than the threshold $T$. In ME, the weights on scales for eigenvalues larger than $T$ are unchanged.\halfblankline

Applying MVI {\reviews with a constant inflation factor} also reduces the contribution of  observation information to the analysis. In contrast to both methods of reconditioning, the reduction in weight is constant for all scales and does not depend on the eigenvectors of $\bfR$. This means that there is no change to the sensitivity to different scales using this method. The analysis will simply pull closer to the background data with the same relative weighting between different observations as occurred for analyses using the original estimated observation error covariance matrix. \halfblankline

We have considered the impact of RR, ME and MVI on the unpreconditioned 3D-Var objective function. For the preconditioned case, 
\citet{johnson05} showed how, when changing the relative weights of the background and observation terms by inflating the ratio of observation and background variances, it is the complex interactions between the error covariance matrices  and the observation operator that affects which scales are present in the analysis. This suggests that in the preconditioned setting MVI will also alter the sensitivity of the analysis to different scales.}

\section{Numerical experiments}\label{sec:NumExpt}

In this section we consider how reconditioning via RR and ME { and application of MVI} affects covariance matrices arising from two different choices of estimated covariance matrices. Both types of covariance matrix are motivated by numerical weather prediction, although similar structures occur for other applications. 

\subsection{Numerical framework}
The first covariance matrix is constructed using a second-order auto-regressive (SOAR) correlation function \citep{Yaglom87} with lengthscale $0.2$ on a unit circle. This correlation function is used in NWP systems \citep{fowler17,stewart13,tabeart17a,waller16,Thiebaux76}
where its long tails approximate the estimated horizontal spatial correlation structure well. In order to construct a SOAR error correlation matrix, $S$, on the finite domain, we follow the method described in \citet{haben11c,tabeart17a}. We consider a one-parameter periodic system on the real line, defined on an equally spaced grid with $N=200$ grid points. We restrict observations to be made only at regularly spaced grid points.
This yields a circulant matrix where the matrix is fully defined by its first row. 
To ensure the corresponding covariance matrix is also circulant, we fix the standard deviation value for all variables 
to be $\sigma = \sqrt{5}$.  \halfblankline

One benefit of using this numerical framework is that it allows us to calculate a simple expression for the update to the standard deviations for ME. We recall that RR updates the variances by a constant, $\delta$. We now show that in the case where $\bfR$ is circulant, ME also updates the variances of $\bfR$ by a constant.\halfblankline

Circulant matrices admit eigenvectors which can be computed directly via a discrete Fourier transform \citep{gray06} (via $\bfR = \bfV\bfLam\bfV^{\dagger}$, where $\dagger$ denotes conjugate transpose). This allows the explicit calculation of the ME standard deviation update given by \eqref{eq:MESigupdate} as
\begin{equation}\label{eq:MEsigCirc}
\bfSig_{ME}(i,i) =\left(\bfSig(i,i)+ \frac{1}{d}\sum_{k=1}^{d}\bfGam(k,k)\right)^{1/2}.
\end{equation}
This follows from \eqref{eq:MESigupdate} because the circulant structure of the SOAR matrix yields $\bfV(i,k)^2 = 1/d$.\halfblankline

We therefore expect both reconditioning methods to increase the SOAR standard deviations by a constant amount. As the original standard deviations were constant, this means that reconditioning will result in constant standard deviations for all variables. These shall be denoted $\sigma_{RR}$ for RR and $\sigma_{ME}$ for ME. { Constant changes to standard deviations also means that an equivalent MVI factor that corresponds to the change can be calculated. This will be denoted by $\alpha$.} \halfblankline

Our second covariance matrix comprises interchannel error correlations for a satellite-based instrument. For this we make use of the Infrared Atmospheric Sounding Instrument (IASI) which is used at many NWP centres within data assimilation systems. A covariance matrix for IASI was diagnosed in 2011 at the Met Office, following the procedure described in \citet{weston11,weston14} (shown in Online Resource 1). The diagnosed matrix was extremely ill-conditioned and required the application of the ridge regression method in order that the correlated covariance matrix could be used in the operational system. We note that we follow the reconditioning procedure of \citet{weston14}, where the reconditioning method is only applied to the subset of $137$ channels that {that are used in the Met Office 4D-Var system. These channels are listed in \citet[Appendix A]{stewart08b}. As the original standard deviation values are not constant across different channels, reconditioning will not change them by a constant amount,  as is the case for Experiment 1. {\reviews We note that the $137\times137$ matrix considered in this paper corresponds to the covariance matrix for one `observation' at a single time and spatial location. The observation error covariance matrix for all observations from this instrument within a single assimilation cycle is a block diagonal matrix, with one block for every observation, each consisting of a submatrix of the $137\times 137$ matrix .} \halfblankline %Therefore, in order to compare reconditioning methods with MVI we use a typical inflation factor $\alpha=1.75$ which is used for this instrument at ECMWF \citep{bormann16}. }\halfblankline%have non-zero off-diagonal correlations. \halfblankline

In the experiments presented in Section \ref{sec:Results} we apply the minimum eigenvalue and the ridge regression methods to both  the SOAR and IASI covariance matrices. The condition number before reconditioning of the SOAR correlation matrix is $81121.71$ and for the IASI matrix we obtain a condition number of $2005.98$. We consider  values of $\kappa_{max}$ in the range $100 - 1000$ for both tests. We note that the equivalence of the minimum eigenvalue method with the minimiser of the Ky Fan $1-d$ norm is satisfied for the SOAR experiment for $\kappa_{max}\ge 168$ and the IASI experiment for $\kappa_{max}\ge 98$.

\subsection{Results}\label{sec:Results}
\subsubsection{Changes to the covariance matrix}
\textbf{Example 1: Horizontal correlations using a SOAR correlation matrix}\\

Due to the specific circulant structure of the SOAR matrix and constant value of standard deviations for all variables, \eqref{eq:RRSigupdate} and \eqref{eq:MEsigCirc} indicate that we expect increases to standard deviations for both methods of reconditioning to be constant. This was found to be the case numerically. In Table~\ref{table:SOARvar} the computed change in standard deviation for different values of $\kappa_{max}$ is given as an absolute value and { as $\alpha$, the multiplicative inflation constant that yields the same change to the standard deviation as each reconditioning method.}
We note that in agreement with the result of Corollary~\ref{lem:RRvsMEvariance} the variance increase is larger for the RR than the ME for all choices of $\kappa_{max}$. Reducing the value of $\kappa_{max}$ increases the  change to standard deviations for both methods of reconditioning.  The increase to standard deviations will result in the observations being down-weighted in the analysis. As this occurs uniformly across all variables for both methods, we expect the analysis to pull closer to the background. Nevertheless, we expect this to be a rather small effect. For this example, even for a small choice of $\kappa_{max}$ the values of the equivalent multiplicative inflation constant, $\alpha$, {\reviews is small, with the largest value of $\alpha = 1.124$ occurring for RR for $\kappa_{max}=100$}.\halfblankline %for both methods are much smaller than the value of $1.75$ that is commonly used (for example at ECMWF \citep{bormann15,bormann16}). }\halfblankline%It is also of interest to consider the corresponding value of the multiplicative inflation constant, $\alpha$. For this example, even for a small choice of $\kappa_{max}$ the values of $\alpha$ for both methods are much smaller than the value of 1.75 that is commonly used (for example at ECMWF \citep{bormann15,bormann16}). The increase to standard deviations will result in the observations being down-weighted in the analysis. As this occurs uniformly across all variables for both methods, we expect the analysis to pull closer to the background.} The updated standard deviations were calculated from the reconditioned SOAR matrix, and agree with theoretical values determined from the update equations \eqref{eq:RRSigupdate} and \eqref{eq:MEsigCirc} to machine precision. \halfblankline
\begin{table}
	\caption{Change in standard deviation for the SOAR covariance matrix for both methods of reconditioning. Columns $4$ and $6$ show $\alpha$, the multiplicative inflation factor corresponding to the values for $\sigma_{RR}$ and $\sigma_{ME}$ respectively.}
	\label{table:SOARvar}
	\begin{tabular}{ l  c  c c c c }
		\hline\noalign{\smallskip}
		$\kappa_{max}$& $\sigma$ & $\sigma_{RR}$ & $\alpha$ corr. RR & $\sigma_{ME}$ & $\alpha$ corr. ME\\
		\noalign{\smallskip}\hline\noalign{\smallskip}
		$1000$&$2.23606 $& $2.26471$&$1.013$& $2.25439$&$1.008$\\
		$500$&$2.23606 $& $2.29340$& $1.026$ & $2.27599$& $1.018$ \\
		$100$&$2.23606 $& $2.51306$& $1.124 $&$2.45737$& $1.099$\\
		\noalign{\smallskip}\hline
	\end{tabular}
\end{table}

As the SOAR matrix is circulant, we can consider the impact of reconditioning on its correlations by focusing on one matrix row. In Figure \ref{fig:SOARcorr} the correlations and percentage change for the 100th row of the SOAR matrix are shown for both methods for $\kappa_{max}=100$.  
These values are calculated directly from the reconditioned matrix. We note that by definition of a correlation matrix,  $\bfC(i,i)=1$ $\forall$ $i$ for all choices of reconditioning. This is the reason for the spike in correlation visible in the centre of Figure \ref{fig:SOARcorr}a and on the right of Figure \ref{fig:SOARcorr}b. { As multiplicative variance inflation does not change the correlation matrix, the black line corresponding to the correlations of the original SOAR matrix also represents the correlations in the case of multiplicative inflation. }We also remark that although ME is not equivalent to the minimiser of the Ky Fan $1-d$ norm for $\kappa_{max}=100$, the qualitative behaviour in terms of correlations and standard deviations is the same for all values in the range $100-1000$. It is important to note that ME is still a well-defined method of reconditioning even if it is not equivalent to the minimiser of the Ky Fan $1-d$ norm. \halfblankline 

\begin{figure}
	\centering
	\subfloat[]{
		\includegraphics[width=75mm]{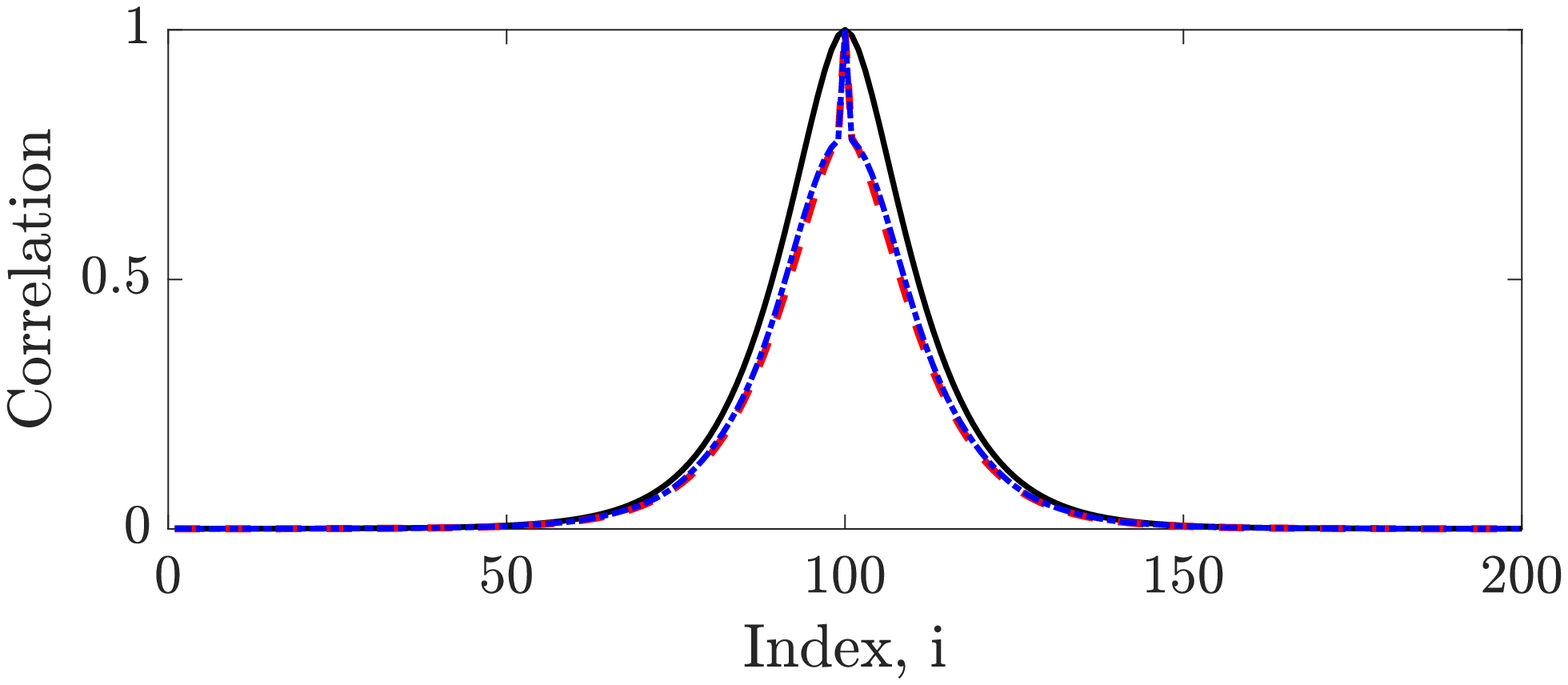}}	
	\subfloat[]{
		\includegraphics[width=75mm]{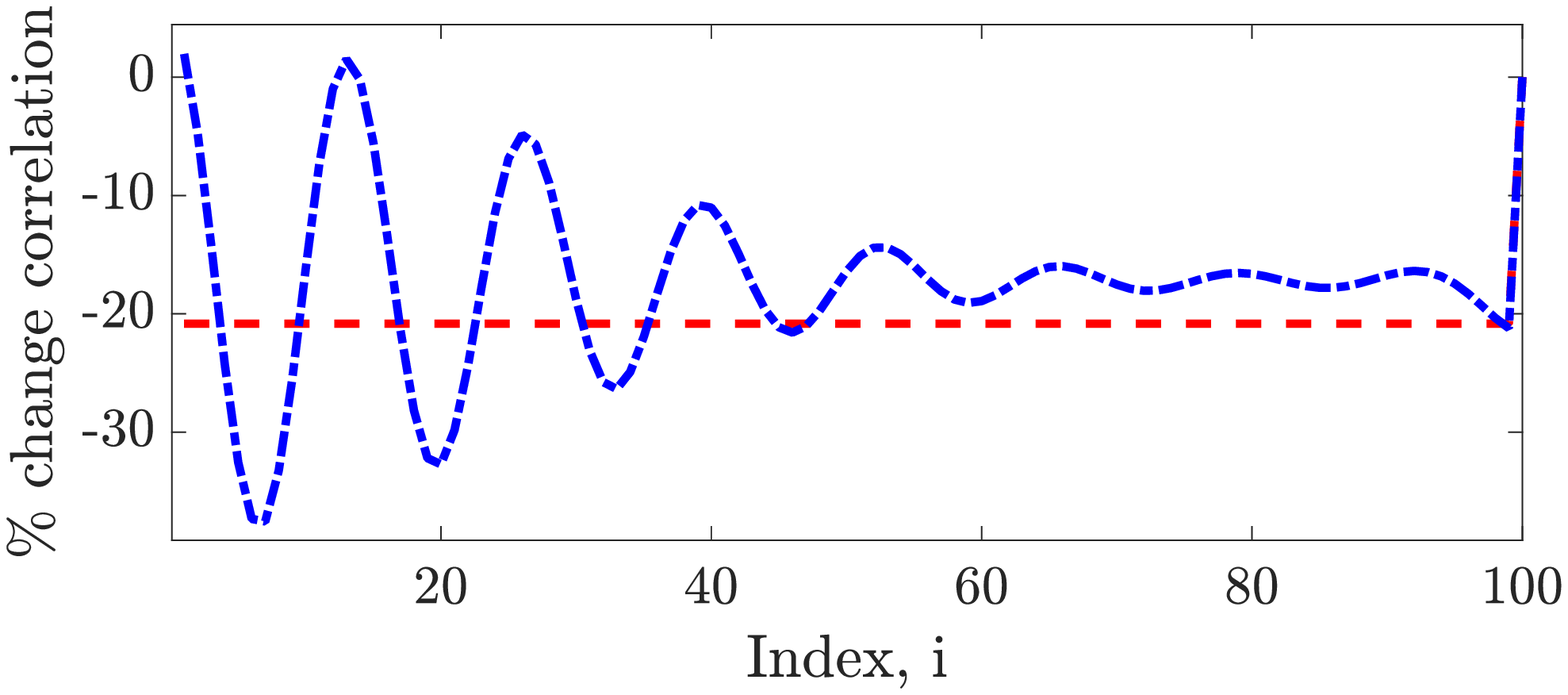}		
	}
	\caption{Changes to correlations between the original SOAR matrix and the reconditioned matrices for $\kappa_{max}=100$. (a) shows $\bfC(100,:)=\bfC_{MVI}(100,:)$ (black solid), $\bfC_{RR}(100,:)$ (red dashed), $\bfC_{ME}(100,:)$ (blue dot-dashed) (b) shows $100\times \frac{\bfC(100,:)-\bfC_{RR}(100,:)}{\bfC(100,:)}$ (red dashed) and $100\times \frac{\bfC(100,:)-\bfC_{ME}(100,:)}{\bfC(100,:)}$ (blue dot-dashed). 
		As the SOAR matrix is symmetric, we only plot the first $100$ entries for (b).} 
	\label{fig:SOARcorr}
\end{figure} 

 Figure \ref{fig:SOARcorr}a shows that for both methods, application of reconditioning reduces the value of off-diagonal correlations for all variables, with the largest absolute reduction occurring for variables closest to the observed variable. Although there is a large change to the off-diagonal correlations, we notice that the correlation lengthscale, which determines the rate of decay of the correlation function, is only reduced by a small amount. This shows that both methods of reconditioning dampen correlation values but do not significantly alter the overall pattern of correlation information.
 Figure \ref{fig:SOARcorr}b shows the percentage change to the original correlation values after reconditioning is applied. 
 For RR, although the difference between the original correlation value and the reconditioned correlation depends on the index $i$, the relative change is constant across all off-diagonal correlations. {As MVI does not alter the correlation matrix, it  would correspond to a horizontal line through $0$ for Figure \ref{fig:SOARcorr}b.}  \halfblankline

  When we directly plot the correlation values for the original and reconditioned matrices in Figure \ref{fig:SOARcorr}a, the change to correlations for ME appears very similar to changes for RR. However, when we consider the percentage change to correlation in Figure \ref{fig:SOARcorr}b we see oscillation in the percentage differences of the ME correlations, showing that the relative effect on some spatially distant variables can be larger than for some spatially close variables. 
  The spatial impact on individual variables differs significantly for this method. We also note that ME increases some correlation values. These are not visible in Figure \ref{fig:SOARcorr} due to entries in the original correlation matrix that are close to zero. Although the differences between $\bfC$ and $\bfC_{ME}$ far from the diagonal are small, small correlation values in the tails of the original SOAR matrix mean that when considering the percentage difference we obtain large values, as seen in Figure \ref{fig:SOARcorr}b.  { This suggests that RR is a more appropriate method to use in this context, as the reconditioned matrix represents the initial correlation function better than ME, where spurious oscillations are introduced. These oscillations occur as ME changes the weighting of eigenvectors of the covariance matrix. As the eigenvectors of circulant matrices can be expressed in terms of discrete Fourier modes, ME has the effect of amplifying the eigenvalues corresponding to the highest frequency eigenvectors. This results in the introduction of spurious oscillations in correlation space. }\halfblankline
  
{ Both methods reduce the correlation lengthscale of the error covariance matrix. In \citet{tabeart17a}, it was shown that reducing the lengthscale of the observation error covariance matrix decreases the condition number of the Hessian of the 3D-Var objective function and results in improved convergence of the minimisation problem. Hence the application of reconditioning methods to the observation error covariance matrix is likely to improve convergence of the overall data assimilation problem. 
	\citet{fowler17} studied the effect on the analysis of complex interactions between the background error correlation lengthscale, the observation error correlation lengthscale and the observation operator in idealised cases. Their findings for a fixed background error covariance, and direct observations, indicate that the effect of reducing the observation error correlation lengthscale (as in the reconditioned cases) is to increase the analysis sensitivity to the observations at larger scales.  In other words, more weight is placed on the large-scale observation information content and less weight on the small scale observation information content. This corresponds with the findings of Section \ref{sec:ObjectiveFn}, where we proved that both methods of reconditioning reduce the weight on small scale observation information in the variational objective function. However, the lengthscale imposed by a more complex observation operator could modify these findings.}\halfblankline %\cite{fowler17} finds that the largest reduction in analysis error (corresponding to the largest conditioning) occurs when observation and background error covariance matrices are sensitive to different scales. As reconditioning reduces the lengthscale of the SOAR matrix, whether reconditioning improves or degrades the analysis depends on the properties of the other error covariance matrix and observation operator. For example if the background correlation lengthscales are larger than observation correlation lengthscales then applying either method of reconditioning to the observation (background) error covariance matrix is likely to improve (degrade) the analysis error.} \halfblankline%Generally, longer correlation lengthscales in observation error covariance matrices are likely to add most benefit \citep{fowler17}. This suggests that applying either method of reconditioning to a SOAR background (observation) error covariance matrix is likely to improve (degrade) the analysis error.}\halfblankline

  \textbf{Example 2: Interchannel correlations using an IASI covariance matrix}

\begin{figure}
	\centering
	\includegraphics[width=84mm]{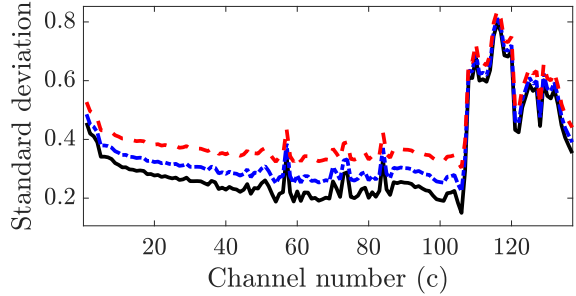}
	\caption{Standard deviations for the IASI covariance matrix $\bfSig$ (black solid), $\bfSig_{RR}$ (red dashed), $\bfSig_{ME}$ (blue dot-dashed) for $\kappa_{max}=100$ and $\alpha = 1.75$. 
		} 
		\label{fig:IASIvar}
\end{figure}

We now consider the impact of reconditioning on the IASI covariance matrix.
We note that there is significant structure in the diagnosed correlations (see \citet[Fig. 8]{stewart14} and online resource 1), with blocks of highly correlated channels in the lower right hand region of the matrix. 
 We now consider how RR, ME and MVI change the variances and correlations of the IASI matrix.\halfblankline

 Figure \ref{fig:IASIvar} shows the standard deviations $\bfSig$, $\bfSig_{RR}$,  $\bfSig_{ME}$ and $\bfSig_{MVI}$. These are calculated from the reconditioned matrices, but the values coincide with the theoretical results of Theorems \ref{lem:updateCandSigRR} and  \ref{lem:changetoSigME}. Standard deviation values for the original diagnosed case have been shown to be close to estimated noise characteristics of the instrument for each of the different channels \citep{stewart14}. 
 { We note that the largest increase to standard deviations occurs for channel 106 only and corresponds to a multiplicative inflation factor for this channel of $2.02$ for RR and $1.81$ for ME.  Channel 106 is sensitive to water vapour and is the channel in the original diagnosed covariance matrix with the smallest standard deviation.
 	 The  choice of $\kappa_{max}=100$ is of a similar size to the value of the parameters used at NWP centres \citep{weston11,weston14,bormann16}. This means that in practice, the contribution of observation information from channels where instrument noise is low is being substantially reduced.\halfblankline} 
   
 Channels are ordered by increasing wavenumber, and are grouped by type. We expect different wavenumbers to have different physical properties, and therefore different covariance structures.
 In particular larger standard deviations are expected for higher wavenumbers due to additional sources of error \citep{weston14}, which is observed on the right hand side of Figure \ref{fig:IASIvar}. For RR, larger increases to standard deviations are seen for channels with smaller standard deviations for the original diagnosed matrix than those with large standard deviations. 
 This also occurs to some extent for ME, although we observe that the update term in \eqref{eq:MESigupdate} is not constant in this case.  {This means that the reduction in weight in the analysis will not be uniform across different channels for ME. } %For channels 1-20, standard deviations are hardly changed by the ME.  
 The result of Corollary \ref{lem:RRvsMEvariance} is satisfied; the increase to the variances is larger for RR than ME. This is particularly evident for channels where the variance from the original diagnosed covariance matrix is small. %  We also see that for all channels bar one, increases to the standard deviations are larger for MVI than for either method of reconditioning. 
{\reviews As MVI increases standard deviations by a constant factor, the largest changes for this method would occur for channels with large standard deviations in the original diagnosed matrix}. This is in contrast to RR, where the largest changes occur for the channels in the original diagnosed matrix with the smallest standard deviation. \halfblankline %We see that even though both methods of reconditioning increase standard deviations, for the majority of channels the increase is much smaller than the equivalent increase for MVI with $\alpha=1.75$. \halfblankline

\begin{figure}
	\centering
	\includegraphics[trim=45mm 0mm 45mm 0mm,width=150mm]{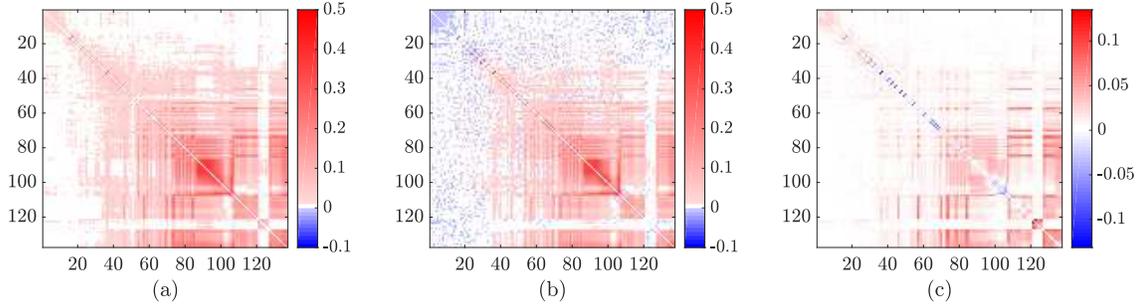}
	\caption{Difference in correlations for IASI (a) $(\bfC - \bfC_{RR})\circ \textit{sign}(\bfC)$, (b) $(\bfC - \bfC_{ME})\circ \textit{sign}(\bfC)$, and  (c) $(\bfC_{ME} - \bfC_{RR})\circ \textit{sign}(\bfC)$, where $\circ$ denotes the Hadamard product. {\reviews Red indicates that the absolute correlation is decreased by reconditioning and blue indicates the absolute correlation is increased. }The colourscale is the same for (a) and (b) but different for (c). Condition numbers of the corresponding covariance matrices are given by $\kappa(\bfR) = 2005.98$, $\kappa(\bfR_{RR}) = 100$ and $\kappa(\bfR_{ME})=100$. } 
	\label{fig:IASIcorr}
\end{figure}

Figure \ref{fig:IASIcorr} shows the difference between the diagnosed correlation matrix, $\bfC$, and the reconditioned correlation matrices $\bfC_{RR}$ and $\bfC_{ME}$. As some correlations in the original IASI matrix are negative, we plot the entries of $(\bfC - \bfC_{RR})\circ \textit{sign}(\bfC)$ and $(\bfC - \bfC_{ME})\circ \textit{sign}(\bfC)$ in Figures \ref{fig:IASIcorr}a and b respectively. Here $\circ$ denotes the Hadamard product, which multiplies matrices of the same dimension elementwise. This allows us to determine whether the magnitude of the correlation value is reduced by the reconditioning method; a positive value indicates that the reconditioning method reduces the magnitude of the correlation, whereas a negative value indicates an increase in the correlation magnitude. For RR, all differences are positive, which agrees with the result of Theorem~\ref{lem:updateCandSigRR}. {As MVI does not change the correlation matrix, an equivalent figure for this method is not given.} We also note that there is a recognisable pattern in Figure \ref{fig:IASIcorr}a, with the largest reductions occurring for the channels in the original diagnosed correlation matrix which were highly correlated. This indicates that this method of reconditioning does not affect all channels equally. \halfblankline

 For ME, we notice that there are a number of entries where the absolute correlations are increased after reconditioning.
 There appears to be some pattern to these entries, with a large number occurring in the upper left hand block of the matrix for channels with the smallest wavenumber \citep{weston14}. However, away from the diagonal for channels 0-40, where changes by RR are very small, the many entries where absolute correlations are increased by ME are much more scattered.  This more noisy change to the correlations could be due to the fact that $96$ eigenvalues are set to be equal to a threshold value by the minimum eigenvalue method in order to attain $\kappa_{max}=100$.
 One method to reduce noise was suggested in \citet{Smith17}, which showed that applying localization methods (typically used to reduce spurious long-distance correlations that arise when using ensemble covariance matrices via the Schur product) after the reconditioning step can act to remove noise while retaining covariance structure. \halfblankline
 
For positive entries, the structure of $\bfC_{ME}$ appears similar to that of RR. 
 	There are some exceptions however, such as the block of channels 121-126  where changes in correlation due to ME are small, but correlations are changed by quite a large amount for RR. The largest elementwise difference between RR and the original diagnosed correlation matrix is $0.138$, whereas the largest elementwise difference between ME and the original diagnosed correlation matrix is $0.0036$. 
 The differences between correlations for ME and RR are shown in Figure 3c. \halfblankline
 
{ For both methods, although the absolute value of all correlations is reduced, correlations for channels 1-70 are eliminated. This has the effect of emphasising the correlations for channels that are sensitive to water vapour. \cite{weston14,bormann16} argue that much of the benefit of introducing correlated observation error for this instrument can be related to the inclusion of correlated error information for water vapour sensitive channels. Therefore, although the changes to the original diagnosed correlation matrix are large it is likely that a lot of the benefit of using correlated observation error matrices is retained. \halfblankline

	We also note that it is more difficult to choose the best reconditioning method in this setting, due to the complex structure of the original diagnosed correlation matrix. {\reviews In particular, improved understanding of how each method alters correlations and standard deviations is not enough to  determine which method will perform best in an assimilation system.} One motivation of reconditioning is to improve convergence of variational data assimilation algorithms. Therefore, one aspect of the system that can be used to select the most appropriate method of reconditioning is the speed of convergence. As ME results in repeated eigenvalues, we would expect faster convergence of conjugate gradient methods applied to the problem $\bfR\bfx = \mathbf{b}$ for $\bfx,\mathbf{b}\in \mathbb{R}^d$ for ME than RR. However,  \citet{campbell16,weston11,weston14,bormann15} find that RR results in faster convergence than ME for operational variational implementations. This is likely due to  interaction between the reconditioned observation error covariance matrix and the observation operator, as the eigenvalues of $\bfH^T\bfR^{-1}\bfH$ are shown to be important for the conditioning of the variational data assimilation problem in \cite{tabeart17a}.\halfblankline
	
	Another aspect of interest is the influence of reconditioning on the analysis and forecast performance. We note that this is likely to be highly system and metric dependent. For example, \citet{campbell16} studies the impact of reconditioning on predictions of meteorological variables (temperature, geopotential height, precipitable water) over lead times from 0 to 5 days. In the U.S. Naval Research Laboratory system, ME performed slightly better at short lead times, whereas RR had improvements at longer lead times \citep{campbell16}. Differences in forecast performance were mixed, whereas convergence was much faster for RR. This meant that the preferred choice was RR. However, in the ECMWF system, \citet{bormann15} studied the standard deviation of first-guess departures against independent observations. Using this metric of analysis skill, ME was found to out-perform RR. The effect of RR on the analysis of the Met Office 1D-Var system is studied in \cite{tabeart17b}, where changes to retrieved variables sensitive to water vapour (humidity, variables sensitive to cloud) are found to be larger than for other meteorological variables such as temperature}.
 %reconditioning eliminates error correlations for channels 1-70. Although all correlations are downweighted, mainly water vapour sensitive channels remain/have zero correlations. Literature: humidity more likely to have correlated observation errors, impact that addition of humidity is likely to have on analysis.
 %True for both methods – more noise chans 1-70 for ME

{\reviews
\subsubsection{Changes to the analysis of a data assimilation problem}
In Section \ref{sec:ObjectiveFn} we considered how the variational objective function is altered by RR, ME and MVI. We found that the two methods of reconditioning reduced the weight on scales corresponding to small eigenvalues by a larger amount than MVI, which changes the weight on all scales uniformly. 
In this section we consider how the analysis of an idealised data assimilation problem is altered by each of the three methods. We also consider how changing $\kappa_{max}$ alters the analysis of the problem. \halfblankline

In order to compare the three methods, we  study how the solution $\bfx$ of a conjugate gradient method applied to the linear system $\bfS\bfx = \mathbf{b}$ changes for RR, ME and MVI, where $\bfS = \bfB^{-1} + \bfH^T\bfR^{-1}\bfH$ is the linearised Hessian associated with the 3D-Var objective function \eqref{eq:CostFn}. To do this we define a `true' solution, $\bfx_{true}$, construct the corresponding data $\mathbf{b}$ and assess how well we are able to recover $\bfx_{true}$ when applying RR, ME and MVI to $\bfS$. 
\citet{haben11c} showed that this is equivalent to solving the 3D-Var objective function in the case of a linear observation operator.
We define a background error covariance matrix, $\bfB\in\mathbb{R}^{200\times 200}$, which is a SOAR correlation matrix on the unit circle with correlation lengthscale $0.2$ and a constant variance of $1$. Our observation operator is given by the identity, meaning that every state variable is observed directly. \halfblankline

 We construct a `true' observation error covariance $\bfR_{true}$, given by a $200$ dimensional SOAR matrix on the unit circle with standard deviation 1 and lengthscale $0.7$. We then sample the Gaussian distribution with zero mean and covariance given by $\bfR_{true}$ 250 times and use these samples to calculate an estimated sample covariance matrix via the Matlab function $cov$. This  estimated sample covariance matrix $\bfR_{est}$ has condition number $\kappa(\bfR_{est})= 3.95\times 10^8$. The largest estimated standard deviation is $1.07$ and the smallest is $0.90$, compared to the true constant standard deviation of $1$.
  RR, ME and MVI are then applied to $\bfR_{est}$ with $\kappa_{max}=100$. When applying MVI, we use two choices of $\alpha$ which correspond to changes to the standard deviations  $(\bfR_{RR}(1,1))^{1/2}$, $\alpha_{RR}= 1.41$, and $(\bfR_{ME}(1,1))^{1/2}$, $\alpha_{ME}=1.39$. The modified error covariance matrices will be denoted $\bfR_{inflRR} = \alpha^2_{RR}\bfR_{est}$ and  $\bfR_{inflME} = \alpha^2_{ME}\bfR_{est}$.  \halfblankline

 We define a true state vector, 
 \begin{equation}\label{eq:xtrue}
\bfx(k)=4\sin(k\pi /100)-5.1\sin(7k\pi /100)+1.5\sin(12k\pi/100)-3\sin(15k\pi/100)+0.75\sin(45k\pi/100),
 \end{equation}
  which has five scales. We then construct $\mathbf{b} \equiv \bfS\bfx$ using $\bfR_{true}$, and apply the Matlab 2018b $pcg.m$ routine to the problem $(\bfB^{-1}+\bfR^{-1})\bfx=\mathbf{b}$ for each choice of $\bfR$. We recall that $\bfS =\bfB^{-1}+\bfH^T\bfR^{-1}\bfH=\bfB^{-1}+\bfR^{-1}$ as $\bfH=\textbf{I}$.  Let $\bfx_{est}$ denote the solution that is found using $\bfR_{est}$ and $\bfx_{mod}$ refer to a solution found using a modified version of $\bfR_{est}$, namely $\bfR_{RR}$, $\bfR_{ME}$, $\bfR_{inflRR}$ or $\bfR_{inflME}$. The maximum number of iterations allowed for the conjugate gradient routine is $200$, and convergence is reached when the relative residual is less than $1\times 10 ^{-6}$.  \halfblankline

\begin{figure}
	\centering

	\includegraphics[trim= 7mm 0mm 15mm 5mm ,clip,width=0.65\textwidth]{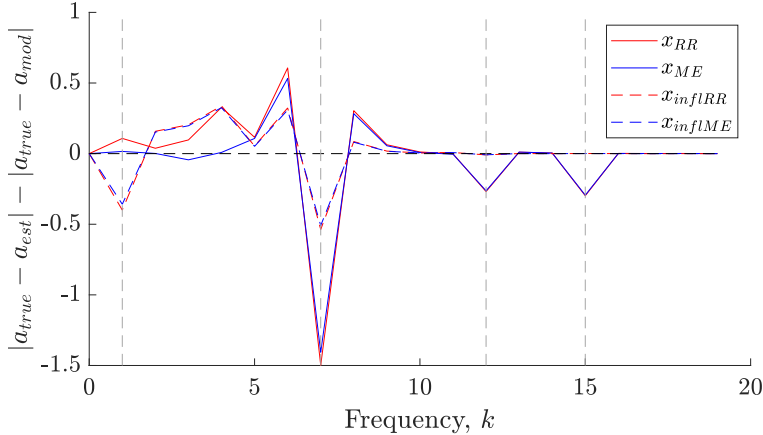}

	\caption{Change in pointwise difference of Discrete Fourier Transform (DFT) from $\bfx_{est}$ to $\bfx_{mod}$ where $a_{est}$ denotes the vector of coefficients of the imaginary part of $DFT(\bfx_{est})$. A positive (negative) value indicates that  $\bfx_{mod}$ is closer to (further from) $\bfx_{true}$ than $\bfx_{est}$ and the amplitude shows how large this change is. Vertical dashed lines show the locations of non-zero values for the true signal.\label{fig:DFT}}
\end{figure}

From Section \ref{sec:ObjectiveFn} we expect RR, ME and MVI to behave differently at small and large scales. We therefore analyse how using each method alters the solution $\bfx$ at different scales using  the discrete Fourier transform (DFT). This allows us to assess how well each scale of $\bfx_{true}$ is recovered for each choice of $\bfR$. As $\bfx_{true}$ is the sum of sine functions, only the imaginary part of the DFT will be non-zero. We therefore define $a_{true}=\text{imag}(DFT(x_{true}))$; similarly $a_{est}=\text{imag}(DFT(x_{est}))$ and $a_{mod}=\text{imag}(DFT(x_{mod}))$. \halfblankline 

By construction, as $\bfx$, given by \eqref{eq:xtrue}, is the sum of sine functions of period $2\pi n/200$ for $n=1,7,12,15,45$, $a_{true}$ returns a signal with 5 peaks, one for each value of $n$ at frequency $k=n$. The amplitude for all other values of $k$ is zero. 
 For frequencies larger than 20, all choices of estimated and modified $\bfR$ recover  $a_{true}$ well. %The largest differences for frequencies 1-20 occur for the second peak where all modified choices of $\bfR$ underestimate the true amplitude, and for frequencies $3-7$, where spurious non-zero amplitudes occur for $\bfx_{est}$ and all choices of $\bfx_{mod}$. 
 Figure \ref{fig:DFT} shows the correction that is applied by the modified choices of $\bfR$ compared to $\bfR_{est}$ for the first 20 frequencies.  A positive (negative) value shows that $a_{mod}$ moves closer to (further from) $a_{true}$ than $a_{est}$. The distance from $0$ shows the size of this change. % from $a_{est}$ to $a_{mod}$. 
 For the first true peak ($k=1$) RR is able to move closer to $a_{true}$ than $a_{est}$. However, both reconditioning methods move further from the truth at the location of true signals $k=7,12,15$. For frequencies where $a_{est}$ has a spurious non-zero signal RR and ME are able to move closer to $a_{true}$ than $a_{est}$. At the location of true signals $k=7,12,15$, MVI makes smaller changes compared to $a_{est}$ than either method of reconditioning. As all modifications to $\bfR_{est}$ move $a_{mod}$ further from $a_{true}$ than $a_{est}$ for $k=7,12,15$, MVI is therefore better able to recover the value of $a_{true}$ than RR or ME at these true peaks. However, MVI introduces a larger error for the first peak at $k=1$ than RR or ME, and changes for frequencies $k>5$ are much smaller than for reconditioning.  
 This agrees with the findings of Section \ref{sec:ObjectiveFn}, that the weight on all scales is changed equally by MVI, whereas both methods of reconditioning result in larger changes to smaller scales and are hence able to make larger changes to amplitudes for higher frequencies. We recall from Section \ref{sec:ObjectiveFn} that ME changes only the smaller scales, whereas RR also makes small changes to the larger scales. This behaviour is seen in Figure \ref{fig:DFT}: for frequencies $k=0$ to $5$ ME results in very small changes, with much larger changes for  frequencies $5\le k\le15$. RR makes larger changes for larger values of $k$, but also moves closer to $a_{true}$ for $1\le k\le3$. \halfblankline

 \begin{figure}
 	\centering
 	\includegraphics[trim= 9mm 0mm 10mm 5mm ,clip,width=0.75\textwidth]{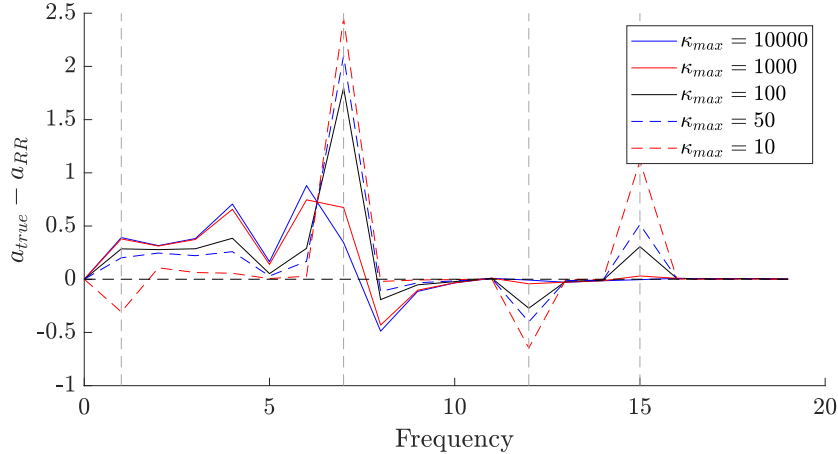}
 	\caption{Difference between $a_{true}$ and $a_{RR}$ for different choices of $\kappa_{max}$. \label{fig:CondChange} Vertical dashed lines show the locations of non-zero values for the true signal.}
 \end{figure}

We now consider how changing $\kappa_{max}$ alters the quality of $\bfx_{RR}$. As the behaviour for $\kappa_{max}=100$ shown in Figure \ref{fig:DFT} was similar for both RR and ME, we only consider changes to RR. Figure \ref{fig:CondChange} shows the difference between $a_{true}$ and $a_{RR}$ for different choices of $\kappa_{max}$. Firstly we consider the true signal that occurs at frequencies $k=1, 7, 12, 15$. For $k=1$ the smallest error occurs for $\kappa_{max}=50$ and the largest error occurs for $\kappa_{max}=10000$. For $k=7,12,15$ the error increases as $\kappa_{max}$ decreases.  For all other frequencies, reducing $\kappa_{max}$  reduces the error in the spurious non-zero amplitudes. For very large values of $\kappa_{max}$ we obtain small errors for the true signal, but larger spurious errors elsewhere. Very small values of $\kappa_{max}$ can control these spurious errors, but fail to recover the correct amplitude for the true signal.  Therefore a larger reconditioning constant will result in larger changes to the analysis. This means that there is a balance to be made in ensuring the true signal is captured, but spurious signal is depressed.  For this framework a choice of $\kappa_{max}=100$ provides a good compromise between recovering the true peaks well and suppressing spurious correlations.%These results agree with the conclusions of Section \ref{sec:ObjectiveFn}, that a larger reconditioning constant will result in larger changes to the objective function. %However, we note that this is a small dimensional problem where we know the true solution exactly. Therefore this analysis may not be possible for more general problems.
 \halfblankline %However, this improvements are much smaller than the degradations that occur to the true amplitudes.  For this framework a choice of $\kappa_{max}=100$ is a good compromise between 

 \begin{table}
	
	\caption{Changes to convergence of RR, MI and MVI for different values of $\kappa_{max}$. For all choices of $\kappa_{max}$, convergence for $\bfR_{true}$ occurs in $17$ iterations and $\bfR_{est}$ occurs in $244$ iterations.}
	\label{table:Convergence}
	\begin{tabular}{ l  c  c c c c }
		\hline\noalign{\smallskip}
		$\kappa_{max}$&  $10,000$ &$1,000$ & $100$ & $50$ & $10$ \\
		\noalign{\smallskip}\hline\noalign{\smallskip}
		RR&$245 $& $244$&$170$& $141$&$73$\\
		ME&$240 $& $239$&$193$& $145$&$76$\\
		Infl RR &$244$&$244$&$238$&$233$&$199$\\
%		Infl ME &$244$&$243$&$238$&$243$&$200$\\
		\noalign{\smallskip}\hline
	\end{tabular}
\end{table}
Finally, Table \ref{table:Convergence} shows how convergence of the conjugate gradient method is altered by the use of reconditioning and MVI. 
Using a larger inflation constant for MVI does lead to slightly faster convergence compared to $\bfR_{est}$.
However, reducing $\kappa_{max}$ leads to a much larger reduction in the number of iterations required for convergence for both RR and ME. %We also see that using ME results in fewer iterations for $\kappa_{max}\ge 1000$, whereas RR results in faster convergence for smaller values of $\kappa_{max}$.   
% However, reconditioning is still required to invert $\bfR$ as $\bfR_{est}$ is rank deficient.
This agrees with results in operational data assimilation systems, where the choice of $\kappa_{max}$ and reconditioning method makes a difference to convergence \cite{weston11,tabeart17b}.

}

\section{Conclusions}\label{sec:Conclusions}

Applications of covariance matrices often arise in high dimensional problems \citep{Pourahmadi13}, such as numerical weather prediction (NWP) \citep{bormann16,weston14}.  In this paper we have examined two methods that are currently used at NWP centres to recondition covariance matrices by altering the spectrum of the original covariance matrix: the ridge regression method, where all eigenvalues are increased by a fixed value, and the minimum eigenvalue method, where eigenvalues smaller than a threshold are increased to equal the threshold value. {We have also considered multiplicative variance inflation, which does not change the condition number or rank of a covariance matrix, but is used at NWP centres \citep{bormann16}.} \halfblankline

 For both reconditioning methods we developed new theory describing how variances are altered. In particular, we showed that both methods will increase variances, and that this increase is larger for the ridge regression method. We also showed that applying the ridge regression method reduces all correlations between different variables. { Comparing the impact of reconditioning methods and multiplicative variance inflation on the variational data assimilation objective function we find that  all methods reduce the weight on observation information in the analysis. However, reconditioning methods have a larger effect on smaller eigenvalues, whereas multiplicative variance inflation does not change the sensitivity of the analysis to different scales.}
We then tested both methods of reconditioning and multiplicative variance inflation numerically on two examples: {Example 1, a spatial covariance matrix, and Example 2, a covariance matrix arising from numerical weather prediction}. In Section \ref{sec:Results} we { illustrated} the theory developed earlier in the work, and also demonstrated that for two contrasting numerical frameworks, the change to the correlations and variances is significantly smaller for the majority of entries for the minimum eigenvalue method. \halfblankline%{\color{red} In Example 1, knowledge of the original correlation function was used to pick the preferred method of reconditioning, whereas in Example 2, the more complex correlation structure meant this was not possible. We also found that standard deviations doubled for some variables for Example 2. It is therefore important for users to understand that the application of reconditioning methods can result in large changes to the correlations and variances of the original covariance matrix. However, for all bar one variable, increases to standard deviations were smaller for both reconditioning methods than multiplicative variance inflation using a typical inflation parameter. }\halfblankline

{ Both reconditioning methods depend on the choice of $\kappa_{max}$, an optimal choice of which will depend on the specific problem in terms of computational resource and required precision. The smaller the choice of $\kappa_{max}$, the more variances and correlations are altered, so it is desirable to select the largest condition number that the  system of interest can deal with. Some aspects of a system that could provide insight into reasonable choices of $\kappa_{max}$ are:
\begin{itemize}
	\item For conjugate gradient methods, the condition number provides an upper bound on the rate of convergence for the problem $A \bfx =\mathbf{b}$ \citep{golub96}, and can provide an indication of the number of iterations required to reach a particular precision \citep{Axelsson96}. Hence $\kappa_{max}$ could be chosen such that a required level of precision is guaranteed for a given number of iterations.
	\item For more general methods, the condition number can provide an indication of the number of digits of accuracy that are lost during computations \citep{Gill86,CheneyWard2013Nmac}.  Knowledge of the error introduced by other system components, {\reviews such as approximations in linearised observation operators and linearised models, relative resolution of the observation network and state variables, precision and calibration of observing instruments,} may give insight into a value of $\kappa_{max}$ that will maintain the level of precision of the overall problem.
	\item The condition number measures how errors in the data are amplified when inverting the matrix of interest \citep{golub96}. Again, the magnitude of errors resulting from other aspects of the system may give an indication of a value of $\kappa_{max}$ that will not dominate the overall precision.
\end{itemize}
%All of these conditions assume that the initial covariance matrix represents the structures of the true covariance matrix, and that the only problem is with conditioning. In the case that some correlations structures are also mis-specified it is likely that other methods, such as inflation, will be required to account for these underestimation of correlation information.
{\reviews For our experiments we considered choices of $\kappa_{max}$ in the range $100-1000$. For Experiment 2 these values are similar to those considered for the same instrument at different NWP centres e.g. 25, 100, 1000 \citep{weston11}, 67 \citep{weston14},  54 and 493 \citep{bormann15}, 169 \citep{campbell16}. We note that the dimension of  this interchannel error covariance matrix 	
	in operational practice is small and only forms a small block of the full observation error covariance matrix.
	   Additionally, the matrix considered in this paper corresponds to one observation type; there are many other observation types with different error characteristics. } \halfblankline 

In this work we have assumed that our estimated covariance matrices represent the desired correlation matrix well, in which case the above conditions on $\kappa_{max}$ can be used. This is not true in general, and it may be that methods such as inflation and localisation are also required in order to constrain the sources of uncertainty that are underestimated or mis-specified. In this case, the guidance we have presented in this paper concerning how to select the most appropriate choice of reconditioning method and target condition number will need to be adapted.  {\reviews Additionally, localisation alters the condition number of a covariance matrix as a side effect; the user does not have the ability to choose the target condition number $\kappa_{max}$ or control changes to the distribution of eigenvalues \citep{Smith17}.} This indicates that reconditioning may still be needed in order to retain valuable correlation information whilst ensuring that the computation of the inverse covariance matrix is feasible. \halfblankline %We showed in Section 3.3 that multiplicative inflation alone is not able to constrain the condition number of a covariance matrix.

The choice of which method is most appropriate for a given situation depends on the system being used and the depth of user knowledge of the characteristics of the error statistics. {\reviews The ridge regression method preserves eigenstructure by increasing the weight of all eigenvalues by the same amount, compared to the minimum eigenvalue method which only increases the weight of small eigenvalues and introduces a large number of repeated eigenvalues.} 
 We have found that ridge regression results in {\reviews constant changes to variances and strict decreases to absolute correlation values,} whereas the minimum eigenvalue method makes smaller, non-monotonic changes to correlations and non-constant changes to variances. }{\reviews In the spatial setting, the minimum eigenvalue method introduced spurious correlations, whereas ridge regression resulted in a constant percentage reduction for all variables. In the inter-channel  case, changes to standard deviations and most correlations were smaller for the minimum eigenvalue method than for ridge regression.}\halfblankline
% This indicates that in situations where users are confident that their estimated covariance matrix represents the true {\reviews correlation structure} well the ridge regression method is better, {\reviews as all correlation values are moved closer to zero}.  The minimum eigenvalue may be more appropriate for situations where standard deviation and correlation values agree with what is expected due to knowledge of the system, but properties of correlation structure (e.g. amplitude, lengthscale, correlation function) are not well understood. \halfblankline} %We also found that although both methods of reconditioning increase standard deviations, for most variables these increases are smaller than changes to standard deviations using multiplicative variance inflation with a commonly used inflation factor. \halfblankline} %EXTRA ECMWF REPORT: For their system, size of sig_o seems to have biggest impact on convergence – but they start with a small condition number.

 Another important property for reconditioning methods is the speed of convergence of minimisation
 	of variational data assimilation problems. {\reviews It is well-known that other aspects of matrix structure, such as repeated or clustered 
 	eigenvalues, are important for the speed of convergence of conjugate gradient minimisation problems. As the condition
 	number is only sensitive to the extreme eigenvalues, conditioning alone cannot fully characterise the
 	expected convergence behaviour. In the data assimilation setting, complex interactions occur between the constituent matrices \citep{tabeart17a}, which can make it hard to determine the best reconditioning method a priori. One example of this is seen for operational implementations in \citet{campbell16,weston11} where the ridge regression method results in fewer iterations for a minimisation procedure than
 	the minimum eigenvalue method, even 
 	though the minimum
 	eigenvalue method yields  observation error covariance matrices with a large number of repeated eigenvalues. 
 	Furthermore, \citet{tabeart17a} found cases in an
 	idealised numerical framework where increasing the condition number of the Hessian  of the data assimilation problem was
 	linked to faster convergence of the minimisation procedure. 
 	Again, this was due to interacting eigenstructures between observation and background terms, which could not be measured by the condition number alone. Additionally, \citet{haben11c,tabeart17a} find that the ratio of background to observation error variance is important for the convergence of a conjugate gradient problem. In the case where observation errors are small, poor performance of conjugate gradient methods is therefore likely.}
 % 	 Complex interactions between background and observation error covariance matrices and the observation operator are described explicitly in \cite{johnson05} for the incremental 4D-Var problem. 
This shows that changes to the analysis of data assimilation problems due to the application of reconditioning methods are likely to be highly system dependent, {\reviews for example due to: quality of estimated covariance matrices, interaction between background and observation error covariance matrices, specific implementations of the assimilation algorithm, and choice of preconditioner and minimisation routine}.  
 	 However, the improved understanding of alterations to correlations and standard deviations for each method of reconditioning provided here may allow users to anticipate changes to the analysis for a particular system of interest using the results from previous idealised and operational studies (e.g. \cite{tabeart17a,fowler17,simonin18,weston14,bormann16}).\halfblankline %{\color{red}Reconditioning methods also alter the problem that is being solved, meaning that the effect on the analysis must also be considered, although that is beyond the scope of this work}.  \halfblankline

% Better understanding and comparison of these two methods is expected to allow users to make informed decisions when selecting methods of reconditioning. This may change based on the desired application; for example the benefits of smaller overall changes to both variances and correlation values observed for the minimum eigenvalue method may be outweighed by large changes to a small number of correlation values, if it is necessary to preserve correlation structure. 

%\begin{appendices}
%	\appendix
	\begin{appendices}
	\section{Equivalence of the minimum eigenvalue method with the Ky Fan 1-d norm method}\label{sec:Append}
\setcounter{equation}{0}
\setcounter{thm}{0}
\numberwithin{equation}{section}
\numberwithin{mydef}{section}
\numberwithin{thm}{section}
We introduce the Ky-Fan $p-k$ norm. We show that the solution to a nearest matrix problem in the Ky-Fan $1-d$ norm, where $\mathbf{X}\in\mathbb{R}^{d\times d}$, is equivalent to the minimum eigenvalue method of reconditioning introduced in 
Section \ref{sec:MinEig} with an additional assumption.
\begin{mydef}
	The Ky Fan $p$-$k$ norm of $\mathbf{X} \in \mathbb{C}^{m\times n}$ is defined as:
	
	\begin{equation}
	||\mathbf{X}||_{p,k} = \left(\sum_{i=1}^{k}\gamma_i(\mathbf{X})^p \right)^{1/p},
	\end{equation}
	where $\gamma_i(\mathbf{X})$ denotes the $i$-th largest singular value of $\mathbf{X}$, $p\ge1$ and $k \in \{1,\dots,\min\{m,n\}\}$.
\end{mydef}
As covariance matrices are positive semi-definite by definition, the singular values of a covariance matrix $\mathbf{X} \in \mathbb{R}^{d\times d}$ are equal to its eigenvalues.
\begin{thm}
	Let $\mathbf{X}\in \mathbb{R}^{d\times d}$ be a symmetric positive semi-definite matrix, with eigenvalues $\lambda_1 \ge \lambda_2 \ge \dots \ge \lambda_d\ge 0$ and corresponding matrix of eigenvectors given by $\mathbf{V}_{\bfR}$. The choice of $\mathbf{\hat{X}}$ that minimises 
	\begin{equation}\label{eq:minprob}
	||\mathbf{X} - \mathbf{\hat{X}}||_{1,p},
	\end{equation}
	subject to the condition	
	$	\kappa(\mathbf{\hat{X}}) = \hat{\kappa},$ for $\hat{\kappa}\ge d-l+1$, 
		is given by
		$	\mathbf{\hat{X}} = \mathbf{V}_{\bfR}Diag(\lambda^*)\mathbf{V}_{\bfR}^T,$
		where $\lambda^*$ is defined by 
	\begin{align}
	\lambda^*_k &= 
	\begin{cases}
	& \mu^* := \frac{\lambda_1}{\hat{\kappa}} \text{if }\lambda_k<\mu^* \\
	& \lambda_k^* = \lambda_k\text{ otherwise.}
	\end{cases} 
	\end{align} 
	and where $l$ is the index such that $\lambda_l\le \mu^* < \lambda_{l-1}$.
\end{thm}
\begin{proof}
	We apply the result given in Theorem 4 of \citep{takana14} for the trace norm (defined as $p=1$ and $k=d$) to find the optimal value of $\mu^*$. Theorem 2 of the same work yields the minimising solution $\mathbf{\hat{X}}$ for the value of $\mu^*$.\halfblankline
	
	We remark that the statement of Theorem 4 of \citep{takana14} uses the stronger assumption that $\hat{\kappa}\ge d$. However, a careful reading of the proof of this theorem indicates that a weaker assumption is sufficient: we assume that $\hat{\kappa}> d-l+1$ where $l$ is the index such that $\lambda_l\le \mu^*<\lambda_{l-1}$.   
\end{proof} 

We note that this optimal value of $\mu^*$ is the same as the threshold $T = \frac{\lambda_1}{\hat{\kappa}}$ defined for the minimum eigenvalue method in \eqref{eq:mineigupdate}
and hence the minimum eigenvalue method is equivalent to the Ky Fan $1$-$d$ minimizer of \eqref{eq:minprob}in the case that $\kappa \ge d-l+1$. \halfblankline

The minimum eigenvalue method is still a valid method of reconditioning when the additional assumption on the eigenvalues of $\mathbf{X}$ is not satisfied. In particular, in the experiments considered in Section \ref{sec:NumExpt} we see qualitatively similar behaviour for the choices of $T$ that satisfy the assumption, and those that do not. It is possible that the lower bound on the condition number imposed by the additional constraint on $\kappa_{max}$ could provide guidance on the selection of the target condition number.

\end{appendices}
\section*{Acknowledgements}
	Thanks to Stefano Migliorini and Fiona Smith at the UK Met Office for the provision of the IASI matrix. This work is funded in part by the EPSRC Centre for Doctoral Training in Mathematics of Planet Earth, the NERC Flooding from Intense Rainfall programme (NE/K008900/1), the EPSRC DARE project (EP/P002331/1) and the NERC National Centre for Earth Observation.

\bibliographystyle{plainnat}
%\bibliographystyle{ac}
%\bibliography{/home/jemima/Dropbox/UoR/PhD/Papers/bibPaper}
\bibliography{bibPaper}

\end{document}

\maketitle
	%\section{Blah}
	\begin{figure}[h]
		\centering
		\includegraphics[width=\textwidth]{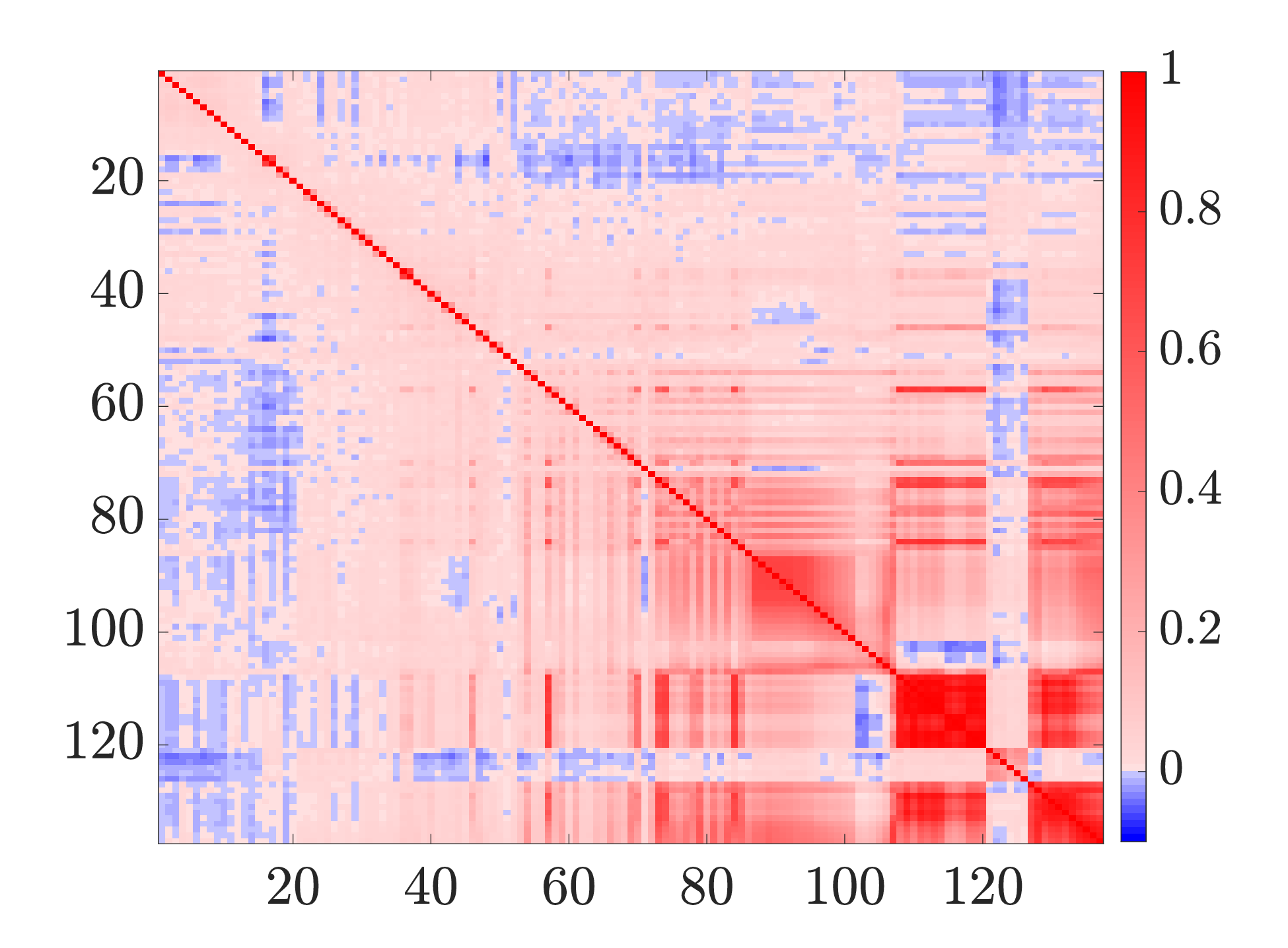}
		\caption{Diagnosed correlation matrix for IASI for the subset of $137$ channels with non-zero off-diagonal entries.}
	\end{figure}